
\input gtmacros
\input gtmonout
\pagenumbers{51}{97}
\volumenumber{1}
\volumeyear{1998}
\volumename{The Epstein birthday schrift}
\received{15 November 1997}
\revised{10 August 1998}
\published{21 October 1998}
\papernumber{3}
\let\qua\stdspace

\def\Endprf{\break\hbox to 0pt{}\endprf}
\def\bksl{\setminus}
\let\empty=\emptyset

\def\text#1{\hbox{\rm #1}}
\let\bbf=\Bbb
\let\script=\cal
\def\brk{\vskip 0pt \noindent}
\def\gap{\vskip 0pt \noindent}

\def\lemma#1{\proclaim{Lemma #1}\rm}
\def\propn#1{\proclaim{Proposition #1}\rm}
\def\theorem#1{\proclaim{Theorem #1}\rm}
\def\thm{\proclaim{Theorem}\rm}
\def\corol#1{\proclaim{Corollary #1}\rm}

\def\defne{\proclaim{Definition}\rm}

\def\em#1{{\it #1\/}}
\def\gap{\vskip 0pt \noindent}

\def\card{\mathop{\rm card}\nolimits}

\def\rhead{\mathop{\rm head}\nolimits}
\def\rtail{\mathop{\rm tail}\nolimits}
\def\diam{\mathop{\rm diam}\nolimits}

\def\dist{\mathop{\rm dist}\nolimits}
\def\Comm{\mathop{\rm Comm}\nolimits}
\def\fr{\mathop{\rm fr}\nolimits}
\def\bet#1{{\langle {#1} \rangle}}

\reflist

\refkey\BeFa
{\bf M Bestvina}, {\bf M Feighn},
{\it Bounding the complexity of simplicial actions on trees},
Invent. Math. 103 (1991) 449--469

\refkey\BeFb
{\bf M Bestvina}, {\bf M Feighn},
{\it A combination theorem for negatively curved\break groups},
J. Differential Geometry 35 (1992) 85--101

\refkey\BeM 
{\bf M Bestvina}, {\bf G Mess},
{\it The boundary of negatively curved groups},
Journal Amer. Math. Soc. 4 (1991) 469--481

\refkey\Boa
{\bf B\,H Bowditch},
{\it Treelike structures arising from continua and convergence groups},
Memoirs Amer. Math. Soc. (to appear)

\refkey\Bob 
{\bf B\,H Bowditch},
{\it Cut points and canonical splittings of hyperbolic groups},
Acta Math. 180 (1998) 145--186

\refkey\Boc 
{\bf B\,H Bowditch},
{\it Group actions on trees and dendrons},
Topology 37 (1998) 1275--1298

\refkey\Bod 
{\bf B\,H Bowditch},
{\it Connectedness properties of limit sets},
Transactions Amer. Math. Soc. (to appear)

\refkey\DeP 
{\bf T Delzant}, {\bf L Potyagailo},
{\it  Accessibilit\'e hi\'erarchique des groupes de
pr\'e-\break sentation
finie},
Strasbourg/Lille preprint (1998)

\refkey\DiD 
{\bf W Dicks}, {\bf M\,J Dunwoody},
{\it Groups acting on graphs},
Cambridge Studies in Advanced Mathematics No. 17,
Cambridge University Press (1989)

\refkey\Du 
{\bf M\,J Dunwoody},
{\it The accessibility of finitely presented groups},
Invent. Math. 81 (1985) 449--457

\refkey\DuS 
{\bf M\,J Dunwoody}, {\bf M\,E Sageev},
{\it JSJ--splittings for finitely presented groups over slender subgroups},
Invent. Math. (to appear)

\refkey\Fe 
{\bf S Ferry},
{\it A stable converse to the Vietoris--Smale theorem with applications
to shape theory},
Transactions Amer. Math. Soc. 261 (1980) 369--386

\refkey \FuP 
{\bf K Fujiwara}, {\bf P Papasoglu},
{\it JSJ decompositions of finitely presented groups and complexes
 of groups}, preprint (1997)

\refkey\GhH 
{\bf E Ghys}, {\bf P de la Harpe},
{\it Sur les groupes hyperboliques d'apr\`es Mikhael Gromov},
Progress in Math. 83, Birkh\"auser (1990)

\refkey\Gr 
{\bf M Gromov},
{\it Hyperbolic groups},
from: ``Essays in Group Theory", S\,M Gersten (editor)
MSRI Publications No. 8, Springer--Verlag (1987) 75--263

\refkey\HoY 
{\bf J\,G Hocking}, {\bf G\,S Young},
{\it Topology},
Addison--Wesley (1961)

\refkey\Ke 
{\bf J\,L Kelley},
{\it General topology},
Graduate Texts in Mathematics 21, Springer--Verlag (reprint of
Van Nostrand edition 1955)

\refkey\Kr 
{\bf J Krasinkiewicz},
{\it Local connectedness and pointed 1--movability},
Bull. Acad. Polon. Sci. S\'er. Sci. Math. Atronom. Phys. 
25 (1977) 1265--1269

\refkey\L 
{\bf G Levitt},
{\it Non-nesting actions on real trees},
Bull. London Math. Soc. 30 (1998) 46--54

\refkey\Ma 
{\bf I Martinez},
{\it Bord d'un produit amalgam\'e sur {\bf Z} de deux groupes libres
 ou de surface},
Orsay preprint (1992)

\refkey\Mi 
{\bf M Mihalik},
{\it Semistability at $\infty $ of finitely generated groups, and
 solvable\break groups},
Topology and its Appl. 24 (1986) 259--264

\refkey\MiT 
{\bf M\,L Mihalik}, {\bf S\,T Tschantz},
{\it Semistability of amalgamated products and\break HNN--extensions},
Memoirs Amer. Math. Soc. No. 471,
Providence, Rhode Island (1992)

\refkey\O 
{\bf J\,P Otal},
{\it Certaines relations d'\'equivalence sur l'ensemble des bouts d'un
groupe libre},
Journal London Math. Soc. 46 (1992) 123--139

\refkey\RS 
{\bf E Rips}, {\bf Z Sela},
{\it Cyclic splittings of finitely presented groups and the canonical
JSJ decomposition},
Annals of Math. 146 (1997) 53--109

\refkey\Se 
{\bf Z Sela},
{\it Structure and rigidity in (Gromov) hyperbolic groups and discrete
 groups in rank $1$ Lie groups II},
Geom. Funct. Anal. 7 (1997) 561--593

\refkey\Sh 
{\bf H Short},
{\it Quasiconvexity and a theorem of Howson's},
from: ``Group theory from a geometrical viewpoint'', E Ghys,
A Haefliger and  A Verjovsky (editors), World Scientific (1991) 168--176

\refkey\St 
{\bf J\,R Stallings},
{\it Group theory and three--dimensional manifolds},
Yale Math. Monographs No. 4, Yale University Press, New Haven (1971)

\refkey\Swa 
{\bf G\,A Swarup},
{\it On the cut point conjecture},
Electron. Res. Announc. Amer. Math. Soc. 2 (1996) 98--100

\refkey\Swe 
{\bf E\,L Swenson},
{\it A cut point tree for a continuum},
preprint (1997)

\refkey\W 
{\bf J\,H\,C Whitehead},
{\it On certain sets of elements in a free group}
Proc. London Math. Soc. 41 (1936) 48--56

\endreflist

\title{Boundaries of strongly accessible\\hyperbolic groups}
\shorttitle{Strongly accessible groups}

\author{B\thinspace H Bowditch}
\asciiauthors{B H Bowditch}
\address{Faculty of Mathematical Studies, University of Southampton\\
Highfield, Southampton SO17 1BJ, Great Britain}
\email{bhb@maths.soton.ac.uk}

\abstract
We consider splittings of groups over finite and two-ended subgroups.
We study the combinatorics of such splittings using generalisations
of Whitehead graphs.
In the case of hyperbolic groups, we relate this to the topology of the
boundary.
In particular, we give a proof that the boundary of a one-ended strongly
accessible hyperbolic group has no global cut point.
\endabstract

\maketitle

\centerline{\small\it 
Dedicated to David Epstein in celebration of his 60th birthday.}

\sectionnumber=-1
\section{Introduction}
\gap\par
In this paper, we consider splittings of groups over finite
and two-ended (ie\ virtually cyclic) groups.
A ``splitting'' of a group, $ \Gamma $, over a class of subgroups may
be viewed a presentation of $ \Gamma $ as a graph of groups, where
each edge group lies in this class.
The splitting is ``non-trivial'' if no vertex group equals $ \Gamma $.
It is said to be a splitting ``relative to'' a given set of subgroups,
if every subgroup in this set can be conjugated into one of the
vertex groups.
Splittings of a given group are often reflected in its large scale
geometry.
Thus, for example, Stallings's theorem [\St] tells us that a finitely
generated group splits non-trivially over a finite group if
and only if it has more than one end.
Furthermore, splittings of a hyperbolic groups over finite and two-ended
subgroups can be seen in the topology of its boundary.
An investigation of this phenomenon will be one of the main objectives of
this paper.
\par
The extent to which a group can be split indefinitely over a certain
class of subgroups is described by the notion of ``accessibility''.
Suppose, $ \Gamma $ is a group, and $ {\script C} $ is a set of subgroups
of $ \Gamma $.
We say that $ \Gamma $ is \em{accessible} over $ {\script C} $ if it
can be represented as a finite graph of groups with all edge groups lying
in $ {\script C} $, and such that no vertex groups splits non-trivially
relative to the incident edge groups.
Dunwoody's theorem [\Du] tells us that any finitely presented group is accessible
over all finite subgroups.
The result of [\BeFa] generalises this to ``small'' subgroups.
\par
There are also stronger notions of accessibility, which have been
considered by Swarup, Dunwoody and others.
One definition is as follows.
Let $ {\script C} $ be a set of subgroups of $ \Gamma $.
Any subgroup of $ \Gamma $ which does not split non-trivially over
$ {\script C} $ is deemed to be ``strongly accessible'' over $ {\script C} $.
Then, inductively, any subgroup which can be expressed as a finite
graph of groups with all edge groups in $ {\script C} $ and all vertex
groups strongly accessible is itself deemed to be ``strongly accessible''.
Put another way, $ \Gamma $ is strongly accessible if some sequence of
splittings of $ \Gamma $ must terminate in a finite number of steps ending up
with a finite number of groups which split no further.
(Of course, this definition leaves open the possibility that there might be a
different sequence of splittings which does not terminate.)
If $ {\script C} $ is the set of finite subgroups, then strong accessibility
coincides with the standard notion of accessibility, and is thus dealt with by
Dunwoody's theorem in the case of finitely presented groups.
Recently Delzant and Potyagailo [\DeP] have shown that any finitely presented
group is strongly accessible over any elementary set of subgroups.
(A set $ {\script C} $ of subgroups is ``elementary'' if no element
of $ {\script C} $ contains a non-cyclic free subgroup, each infinite
element of $ {\script C} $ is contained in a unique maximal element of
$ {\script C} $, and each maximal element of $ {\script C} $ is equal to
its normaliser in $ \Gamma $.)
\par
If $ \Gamma $ is hyperbolic in the sense of Gromov [\Gr], then the set of
all finite and two-ended subgroups is elementary.
Thus, the result of [\DeP] tells us that $ \Gamma $ is strongly accessible.
(In the context of hyperbolic groups, we shall always take ``strongly
accessible'' to mean strongly accessible over finite and two-ended subgroups.)
\par
The boundary, $ \partial \Gamma $, of $ \Gamma $ is a compact metrisable
space, and is connected if and only if $ \Gamma $ is one-ended.
In this case, it was shown in [\BeM] that $ \partial \Gamma $ is locally
connected provided it has no global cut point.
In this paper, we show (Theorem 9.3):
\gap
\thm
{\sl
The boundary of a one-ended strongly accessible group has no global cut
point.
}
\gap
\ppar
Thus, together with [\DeP] and [\BeM], we arrive at the conclusion that
the boundary of every one-ended hyperbolic group is locally connected.
This was already obtained by Swarup [\Swa] using results from
[\Boa,\Boc,\L] shortly after the original draft of this paper was circulated
(and prior to the result of [\DeP]).
An elaboration of the argument was given shortly afterwards in [\Bod].
\par
One consequence of this local connectedness is the fact that
every hyperbolic group is semistable at infinity [\Mi].
(It has been conjectured that every finitely presented group has this
property.)
This implication was observed by Geoghegen and reported in [\BeM].
I am indebted to Ross Geoghegen for the following elaboration of how this
works.
The semistability of an accessible group is equivalent to
the semistability of each of its maximal one-ended subgroups.
Suppose, then, that $ \Gamma $ is a one-ended hyperbolic group.
It was shown in [\BeM] that $ \partial \Gamma $ naturally compactifies
the Rips complex, so as to give a contractable ANR, with
$ \partial \Gamma $ embedded as a Z--set.
It follows that semistability at infinity for $ \Gamma $ is
equivalent to $ \partial \Gamma $ being pointed 1--movable, the
latter property being intrinsic to $ \partial \Gamma $.
Moreover, it was shown in [\Kr] that a metrisable continuum is pointed
1--movable if and only if it has the shape of a Peano continuum
(see also [\Fe]).
It follows that if $ \Gamma $ is one-ended hyperbolic, then
$ \partial \Gamma $ is semistable at infinity if and only of
$ \partial \Gamma $ has the shape of a Peano continuum.
(We remark that an alternative route to semistability for a hyperbolic
group would be to use the result of [\MiT] in place of Theorem 8.1 of
this paper, together with the results of [\Boa,\Boc].)
\par
We shall carry out much of our analysis of splitting in a fairly general
context.
We remark that any one-ended finitely presented group admits a canonical
splitting over two ended subgroups, namely the JSJ splitting
(see [\RS,\DuS,\FuP], or in the context of hyperbolic groups [\Se,\Bob]).
The vertex group are again finitely presented, and so we can split them over
finite subgroups as necessary and iterate the process, discarding any finite
vertex groups that arise along the way.
This eventually leads to a canonical decomposition of the group into
one-ended subgroups, none of which split over any two-ended subgroup.
Further discussion of this procedure will be given in Section 9.
We shall not make any explicit use of the JSJ splitting in this paper.
\par
In this paper, we shall be considering in some detail the general issue of
splittings over two-ended subgroups.
One point to note (Theorem 2.3) is the following:
\gap
\thm
{\sl
The fundamental group of a finite graph of groups with two-ended
edge groups is one-ended if and only if no vertex group splits over
a finite subgroup relative to the incident edge groups.
}
\gap\ppar
(The case where the vertex groups are all free or surface groups is
dealt with in [\Ma].)
\par
To find a criterion for recognising whether a given group splits over
a finite group relative to a given finite set of two-ended subgroups, we
shall generalise work of Whitehead and Otal in the case of free groups.
Given a free group, $ F $, and a non-trivial element, $ \gamma \in F $,
we say that $ \gamma $ is ``indecomposable'' in $ F $, if it cannot
be conjugated into any proper free factor of $ F $.
\par
This can be interpreted topologically.
Note that the boundary, $ \partial F $, of $ F $ is a Cantor set.
We define an equivalence relation, $ \mathord{\approx} $, on
$ \partial F $, by deeming that $ x \approx y $ if and only if
either $ x=y $ or $ x $ and $ y $ are the fixed points of some
conjugate of $ \gamma $.
Now, it's easily verified that this relation is closed, and so the
(equivariant) quotient, $ \partial F/\mathord{\approx} $ is compact
hausdorff.
It was shown in [\O] that $ \gamma $ is indecomposable if and only
if $ \partial F/\mathord{\approx} $ is connected
(in which case, $ \partial F/\mathord{\approx} $
is locally connected and has no global cut point).
\par
A combinatorial criterion for indecomposability is formulated in [\W].
Let $ a_1, a_2 ,$ $\ldots, a_n $ be a system of free generators for $ F $.
Let $ w $ be a reduced cyclic word in the $ a_i $'s and their inverses
representing (the conjugacy class of) $ \gamma $.
Let $ {\script G} $ be the graph (called the ``Whitehead graph'') with
vertex set $ a_1 ,\ldots, a_n, a_1^{-1} ,\ldots, a_n^{-1} $, and with
$ a_i^{\epsilon_i} $ deemed to be adjacent to $ a_j^{\epsilon_j} $
if and only if the string $ a_i^{\epsilon_i} a_j^{-\epsilon_j} $
occurs somewhere in $ w $ (where $ \epsilon_i, \epsilon_j \in \{ -1,1 \} $).
Suppose we choose the generating set so as to minimise the length
of the word $ w $.
Then (a simple consequence of) Whitehead's lemma tells us that
$ \gamma $ is indecomposable if and only if $ {\script G} $ is connected.
(Moreover in such a case, $ {\script G} $ has no cut vertex.)
\par
This can be reinterpreted in terms of what we shall call ``arc systems''.
Let $ T $ be the Cayley graph of $ F $ with respect to free generators
$ a_1 \ldots a_n $.
Thus, $ T $ is a simplicial tree, whose ideal boundary, $ \partial T $, may
be naturally identified with $ \partial F $.
The element $ \gamma $ determines a biinfinite arc, $ \beta $, in $ T $,
namely the axis of $ \gamma $.
Let $ {\script B} $ be the set of images of $ \beta $ under $ \Gamma $.
We refer to $ {\script B} $ as a ($ \Gamma $--invariant) ``arc system''.
We can reconstruct the Whitehead graph, as well as the equivalence relation
$ \mathord{\approx} $, from this arc system in a simple
combinatorial fashion, as described in Section 3.
The above discussion applies equally well if we
replace $ \gamma $ by a finite set, $ \{ \gamma_1 ,\ldots, \gamma_p \} $,
of non-trivial elements of $ \Gamma $.
\par
One can generalise these notions to an arbitrary hyperbolic group,
$ \Gamma $.
Suppose that $ \{ H_1 ,\ldots, H_p \} $ is a finite set of two-ended
subgroups of $ \Gamma $.
We define an equivalence relation, $ \mathord{\approx} $, on
$ \partial \Gamma $ by identifying the two endpoints of each
conjugate to each $ H_i $.
Thus, as before, $ \partial \Gamma/\mathord{\approx} $ is hausdorff.
We shall see (Theorem 5.2) that:
\gap
\thm
{\sl
$ \partial \Gamma/\mathord{\approx} $ is connected if and only if
$ \Gamma $ does not split over a finite group relative to
$ \{ H_1 ,\ldots, H_p \} $.
}
\gap
\ppar
We can also give a combinatorial means of recognising if $ \Gamma $
splits in this way.
We can decompose its boundary, $ \partial \Gamma $, as a disjoint
union of two $ \Gamma $--invariant sets, $ \partial_0 \Gamma $ and
$ \partial_\infty \Gamma $, where $ \partial_\infty \Gamma $ is the
set of singleton components of $ \partial \Gamma $.
Algebraically this corresponds the action of $ \Gamma $ on a simplicial
tree, $ T $, without edge inversions, with finite quotient, and with finite
edge stabilisers and finite or one-ended vertex stabilisers.
Such an action is given by the accessibility theorem [\Du].
Each of the vertex groups is quasiconvex, and hence intrinsically
hyperbolic.
Now, $ \partial_\infty \Gamma $ can be canonically identified with
$ \partial T $, and the connected components of $ \partial_0 \Gamma $
are precisely the boundaries of the infinite vertex stabilisers.
The infinite vertex stabilisers are, in fact, precisely the maximal
one-ended subgroups of $ \Gamma $.
(Note that $ \Gamma $ is virtually free if and only of
$ \partial_0 \Gamma = \emptyset $.)
We can construct an analogue of the Whitehead graph by considering
the arc system on $ T $, consisting of all the translates of the
axes of those $ H_i $ which do not fix any vertex of $ T $.
\par
This combinatorial construction can be carried out for
any group which is accessible over finite subgroups.
Put together with Theorem 2.3, this gives a combinatorial criterion
for recognising when a finitely presented group represented as graph
of groups with two-ended edge groups is one-ended.
This generalises work of Martinez [\Ma].
It is also worth remarking that the result of [\BeFb] tells us that
such a group is hyperbolic if and only if all the vertex groups are
hyperbolic, and there is no Baumslag--Solitar (or free abelian) subgroup.
\par
The structure of this paper is roughly as follows.
In Section 1, we explore some general facts about groups accessible over
finite groups.
In Section 2, we give a criterion (Theorem 2.3) for a finite graph of
groups with two-ended edge groups to be one-ended.
In Section 3, we study arc systems on trees and their connections to
Whitehead graphs.
In Section 4, we give an overview of some general facts about quasiconvex
splittings.
In Section 5, we look at certain quotients of the boundaries of hyperbolic
groups, and relate this to some of the combinatorial results of Section 3.
In Section 6, we set up some of the general machinery for analysing
the topology of the boundaries of hyperbolic groups which split over
two-ended subgroups.
In Section 7, we look at some implications concerning connectedness
properties of boundaries.
In Section 8, we apply this specifically to global cut points.
Finally, in Section 9, we discuss further the question of strong accessibility
of groups over finite and two-ended subgroups.
\par
Much of the material of the original version of this paper was worked out
while visiting the University of Auckland.
The first draft was written at the University of Melbourne.
I would like to thank Gaven Martin as well as Craig Hodgson and
Walter Neumann for their respective invitations.
The paper was substantially revised in Southampton, with much of the material
of Sections 1, 2, 3 and 5 added.
I am also grateful to Martin Dunwoody for helpful conversations regarding the
latter.
Ultimately, as always, I am indebted to my ex-PhD supervisor David Epstein
for first introducing me to matters hyperbolical.
\gap\gap
\section{Trees and splittings}
\gap\par
In this section, we introduce some terminology and notation relating
to simplicial trees and group splittings.
\par
Let $ T $ be a simplicial tree, which we regard a 1--dimensional CW--complex.
We write $ V(T) $ and $ E(T) $ respectively for the vertex set and
edge set.
Given $ v,w \in V(T) $, we write $ \dist(v,w) $ for the distance between
$ v $ and $ w $, in other words, the number of edges in the arc connecting
$ v $ to $ w $.
If $ {\vec e} \in {\vec E}(T) $ and $ v \in V(T) $, we say that
$ {\vec e} $ ``points towards'' $ v $ if $ \dist(v,\rtail({\vec e}))
= \dist(v,\rhead({\vec e}))+1 $.
\par
If $ S \subseteq T $ is a subgraph, we write $ V(S) \subseteq V(T) $ and
$ E(S) \subseteq E(T) $ for the corresponding vertex and edge sets.
A subtree of $ T $ is a connected subgraph.
Of particular interest are ``rays'' and ``biinfinite arcs''
(properly embedded subsets homeomorphic to $ [0,\infty) $ and
$ {\bbf R} $ respectively.)
\par
We may define the ideal boundary, $ \partial T $, of $ T $, as the
set of cofinality classes of rays in $ \Sigma $.
We shall only be interested in $ \partial T $ as a set.
(In fact, $ T \cup \partial T $ can be given a natural compact topology as a
dendron, as discussed in [\Boa].
It can also be given a finer topology by viewing $ T $ has a Gromov
hyperbolic space, and $ \partial T $ as its Gromov boundary.)
If $ S \subseteq T $ is a subgraph, we write $ \partial S \subseteq
\partial T $ for the subset arising from those rays which lie in
$ S $.
Note that if $ \beta $ is a biinfinite arc, then $ \partial \beta $
contains precisely two points, $ x,y \in \partial T $.
We say that $ \beta $ {\it connects\/} $ x $ to $ y $.
\par
Further discussion of general simplicial trees will be given in
Sections 2 and 3.
We now move on to consider group actions on trees.
\par
Let $ G $ be a group.
A $ G $--\em{tree} is a simplicial tree, $ T $, admitting a simplicial
action of $ G $ without edge inversions.
If $ v \in V(T) $ and $ e \in E(T) $, we write $ G_T(v) $ and
$ G_T(e) $ for the corresponding vertex and edge stabilisers respectively.
Where there can be no confusion, we shall abbreviate these to
$ G(v) $ and $ G(e) $.
Such a tree gives rise to a splitting of $ G $ as a graph of groups,
$ G/T $.
We shall say that $ T $ is \em{cofinite} if $ T/G $ is finite.
We shall usually assume that $ T $ is \em{minimal}, ie\ that there is
no proper $ G $--invariant subtree.
This is the same as saying that $ T $ has no terminal vertex, or, on
the level of the splitting, that no vertex group of degree one is equal to
the incident edge groups.
Such a vertex will be referred to as a \em{trivial vertex}.
A subset (usually a subgroup) $ H $, of $ G $ is \em{elliptic} with
respect to $ T $, if it lies inside some vertex stabiliser.
If $ {\script H} $ is a set of subsets of $ G $, we say that the splitting
is \em{relative to} $ {\script H} $, if every element of $ {\script H} $ is an
elliptic subset.
We note that any finite subgroup of a group is elliptic with respect to
every splitting.
Thus any splitting of any group is necesarily relative to the set of
all finite subgroups.
\par
Suppose that $ F $ is a $ G $--invariant subgraph of $ T $, we can
obtain a new $ G $--tree, $ \Sigma $, by collapsing each component of
$ F $ to a point.
We speak of the splitting $ T/G $ as being a \em{refinement} of the
splitting $ \Sigma/G $.
Note that one may obtain a refinement of a given graph of groups, if
one of the vertex groups splits relative to its incident edge groups.
\par
We say that a $ G $--tree, $ T' $, is a \em{subdivision} of $ T $, if
it is obtained by inserting degree--2 vertices into the edges of $ T $ in
a $ G $--equivariant fashion.
Suppose that $ \Sigma $ is another $ G $--tree.
A \em{folding} of $ T $ onto $ \Sigma $ is a $ G $--equivariant map
of $ T $ onto $ \Sigma $ such that each edge of $ T $ either gets mapped
homeomorphically onto an edge of $ \Sigma $ or gets collapsed to a vertex
of $ \Sigma $.
A \em{morphism} of $ T $ onto $ \Sigma $ is a folding of some subdivision
of $ T $.
Such maps are necessarily surjective provided that $ \Sigma $ is minimal.
Clearly a composition of morphisms is a morphism.
\par
We say that $ T $ \em{dominates} $ \Sigma $ (or that the splitting
$ T/G $ \em{dominates} $ \Sigma/G $) if there exists a morphism from
$ T $ to $ \Sigma $.
It's not hard to see that this is equivalent to saying that every
vertex stabiliser in $ T $ is elliptic with respect to $ \Sigma $.
We say that $ T $ and $ \Sigma $ are \em{equivalent} if each dominates
the other.
This is equivalent to saying that a subset of $ G $ is elliptic with
respect to $ T $ if and only if it is elliptic with respect to $ \Sigma $.
\par
Suppose that $ T $ is cofinite.
If $ T $ dominates $ \Sigma $, then $ \Sigma $ is also cofinite.
In this case, any morphism from $ T $ to $ \Sigma $ expands combinatorial
distances by at most a bounded factor (namely the maximum number of edges into
which we need to subdivide a given edge of $ T $ to get a folding.)
Also, any two morphisms remain a bounded distance apart.
In particular, any self-morphism of a cofinite tree is a bounded distance
from the identity map, and is thus a quasiisometry.
Suppose that $ T $ and $ \Sigma $ are equivalent, and that
$ \phi\co T\longrightarrow \Sigma $ is a morphism.
Let $ \psi\co \Sigma \longrightarrow T $ be any morphism.
Now, since $ \psi $ expands distances by a bounded factor, and
$ \psi \circ \phi $ is a quasiisometry, it follows that $ \phi $ is itself
a quasiisometry.
In summary, we have shown:
\gap
\lemma{1.1}
{\sl
If $ T $ and $ \Sigma $ are equivalent cofinite $ G $--trees, then
any morphism from $ T $ to $ \Sigma $ is quasiisometry.
}
\endprf
\par
We see from the above discussion that there is a natural bijective
correspondence between the boundaries, $ \partial T $ and $ \partial \Sigma $,
of $ T $ and $ \Sigma $.
\gap
\lemma{1.2}
{\sl
Suppose that $ T $ and $ \Sigma $ are cofinite $ G $--trees with
finite edge-stabilisers.
If $ \phi\co T\longrightarrow \Sigma $ is a folding, then only finitely
many edges of $ T $ get mapped homeomorhically under $ \phi $ to
any given edge of $ \Sigma $.
}
\gap
\prf
If $ \gamma \in \Gamma $ and $ e, \gamma e \in E(T) $ both get mapped
homeomorpically onto some edge $ \epsilon \in E(\Sigma) $, then
$ \gamma \in \Gamma_\Sigma(\epsilon) $.
There are thus only finitely many such edges in the $ \Gamma $--orbit
of $ e $ in $ E(T) $.
The result follows since $ E(T)/\Gamma $ is finite.\Endprf
\par
We shall need to elaborate a little on the notion of accessibility over
finite groups.
For the remainder of this section, all splittings will be assumed
to be over finite groups, and the term ``accessible'' is assumed to mean
``accessible over finite groups''.
\par
We shall say that a graph of groups is \em{reduced} if no vertex group
of degree one or two is equal to an incident edge group.
(Every graph of groups is a refinement of a reduced graph.)
We say that a group $ G $ is ``accessible'' if there is a bound on the
complexity (as measured by the number of edges) of a splitting of
$ G $ as a reduced graph of groups (with finite edge groups).
Among graphs of maximal complexity, one for which the sum of the
orders of the edge stabilisers is minimal will be referred to as
a ``complete splitting''.
By Dunwoody's theorem [\Du], any finitely presented group is accessible.
(This has been generalised to splittings over small subgroups by
Bestvina and Feighn [\BeFa].)
\par
This can be rephrased in terms of one-ended subgroups.
For this purpose, we define a group to be \em{one-ended} if it is infinite and
does not split non-trivially (over any finite subgroup).
Thus, by Stallings's theorem, this coincides with the usual topological
notion for finitely generated groups.
Suppose that $ G $ is accessible, and we take a complete splitting of $ G $.
Now any splitting of a vertex group is necessarily relative to the incident
edge groups, and so would give rise to a refined splitting.
It is possible that this refined splitting may no longer be reduced, but
in such a case, we can coalesce two vertex groups, to produce a reduced graph
with one smaller edge stabiliser than the original, thereby contradicting
completeness.
In summary, we see that all the vertex groups of a complete splitting are
either finite or one-ended.
In fact, we see that the infinite vertex groups are precisely the
maximal one-ended subgroups.
It turns out that there is a converse to this statement: any group which
can be represented as a finite graph of groups with finite edge groups
and with all vertex groups finite or one-ended is necessarily accessible
(see [\DiD]).
\par
Finally, suppose that $ G $ is accessible, and we represent it as
a finite graph of groups over finite subgroups.
Now each vertex group must be accessible.
Taking complete splittings of each of the vertex groups, we can see
that we can refine the original splitting in such a way that all
the vertex groups are finite or one-ended.
(It is possible that this refinement might not be reduced.)
\par
Now, let $ G $ be an accessible group, and let $ T $ be a cofinite tree
with finite edge stabilisers and with every vertex stabilisers
either finite or one-ended.
The infinite vertex groups are canonically determined.
We have also observed that finite groups are always elliptic in
any splitting.
It follows that if $ T' $ is another such $ G $--tree, then
$ T $ and $ T' $ are equivalent, by Lemma 1.1.
In particular $ \partial T $ and $ \partial T' $ can be canonically
(and hence $ G $--equivariantly) identified.
We can thus associate to any accessible group, $ G $, a canonical
$ G $--set, $ \partial_\infty G $, which we may identify with the boundary
of any such $ G $--tree.
\par
Clearly in the case of a free group, we just recover the usual boundary.
More generally, if $ G $ is (word) hyperbolic (and hence accessible)
then we may identify $ \partial_\infty G $ with the set of singleton
components of the boundary, $ \partial G $.
In fact, as discussed in the introduction, we can write
$ \partial G $ as a disjoint union $ \partial_0 G \sqcup \partial_\infty G $,
where each component of $ \partial_0 G $ is the boundary of a maximal
one-ended subgroup of $ G $.
\par
We shall make some further observations about accessible groups in
connection with strong accessibility in Section 9.
\gap\gap
\section{Splittings over two-ended subgroups}
\gap\par
The main aim of this section will be to give a proof of Theorem 2.3.
We first introduce some terminology regarding ``arc systems'' which will
be relevant to later sections.
\par
Let $ T $ be a simplicial tree.
\gap
\defne
An \em{arc system}, $ {\script B} $, on $ T $ consists of a set of
biinfinite arcs in $ T $.
\brk
We say that $ {\script B} $ is \em{edge-finite} if at most finitely
many elements of $ {\script B} $ contain any given edge of $ T $.
\brk
If $ G $ is a group, and $ T $ is a $ G $--tree, then we shall assume that
an arc system on $ T $ is $ G $--invariant.
\gap
\par
Recall that a subgroup, $ H $, of $ G $ is ``elliptic'' if it fixes a
vertex of $ T $.
If $ H $ is two-ended (ie\ virtually cyclic) then either $ H $ is
elliptic, or else there is a biinfinite $ \beta $ in $ T $ which is
$ H $--invariant.
In the latter case, we say that $ H $ is \em{hyperbolic} and that
$ \beta $ is the \em{axis} of $ H $.
Clearly, the $ H $--stabiliser of any edge of $ {\script B} $ is finite.
\par
Suppose now that all edge stabilisers of $ T $ are finite.
Then every hyperbolic two-ended subgroup of $ G $ lies in a unique maximal
two-ended subgroup of $ G $, namely the setwise stabiliser of the axis.
Note also that there are only finitely many two-ended subgroups, $ H $, with
a given axis, $ {\script B} $, and with the number of edges of
$ \beta/H $ bounded.
In particular, we see that only finitely many $ G $--conjugates of
a given hyperbolic two-ended subgroup, $ H $, can share the same axis.
\par
Suppose, now, that $ {\script H} $ is a finite union of conjugacy
classes of two-ended subgroups of $ G $, and that $ {\script B} $ is
the set of all axes of all hyperbolic elements of $ {\script H} $.
(In other words, $ {\script B} $ is an arc-system with $ {\script B}/\Gamma $
finite, and such that the setwise stabiliser of each element of
$ {\script B} $ is infinite, and hence two-ended.)
We note:
\gap
\lemma{2.1}
{\sl
The arc system $ {\script B} $ is edge-finite.
}
\gap
\prf
We want to show that any given edge lies in a finite number of
elements of $ {\script B} $.
Without loss of generality, we can suppose that $ {\script B} $ conists
of the orbit of a single arc, $ \beta $.
Let $ H $ be the setwise stabiliser of $ \beta $.
Choose any edge $ e \in T $.
Let $ K \le G $ be the stabiliser of $ e $.
Without loss of generality, we may as well suppose that
$ e \in E(\beta) $.
Note that $ E(\beta)/H $ is finite.
Now, the $ G $--orbit, $ Ge $, of $ e $ meets $ E(\beta) $ in an
$ H $--invariant set consisting of finitely many $ H $--orbits,
say $ Ge \cap E(\beta) = Hg_1 e \cup Hg_2 e \cup \cdots \cup Hg_n e $, where
$ g_i \in G $.
\par
Suppose that $ e \subseteq g \beta $, for some $ g \in G $.
Now $ g^{-1} e \in E(\beta) $, so $ g^{-1} e = h g_i e $ for some
$ h \in H $, and $ i \in \{ 1 ,\ldots, n \} $.
Thus $ ghg_i \in K $, so $ gH = kg_i^{-1} H $ for some $ k \in K $.
Since $ K $ is finite, there are finitely many possibilities for the
right coset $ gH $, and hence for the arc $ g \beta $.
\endprf
\par
Now, let $ {\script H} $ be any finite union of conjugacy classes of
two ended subgroups of $ G $, as above.
Recall that to say that $ G $ splits over a finite subgroup relative to
$ {\script H} $ means that there is a non-trivial $ G $--tree with finite edge
stabilisers, and with each element of $ {\script H} $ elliptic with
respect to $ T $.
We can always take such a $ G $--tree to be cofinite, and indeed to have
only one orbit of edges.
We say that $ {\script H} $ is \em{indecomposable} if $ G $ does not split
over any finite group relative to $ {\script H} $.
\par
In Section 3, we shall give a general criterion for indecomposability in
terms of arc systems.
For the moment, we note:
\gap
\lemma{2.2}
{\sl
Suppose that $ G $ is a group and that $ T $ is a $ G $--tree with finite
edge stabilisers.
Suppose that $ {\script H} $ is a finite union of conjugacy classes of
two-ended subgroups of $ G $.
Let $ {\script B} $ be the arc system consisting of the set of axes of
hyperbolic elements of $ G $.
If $ {\script H} $ is indecomposable, then each edge of $ T $ lies in at
least two elements of $ {\script B} $.
}
\gap
\prf
Suppose that $ T \ne \bigcup {\script B} $.
Then, collapsing each component of $ \bigcup {\script B} $ to a point,
we obtain another $ G $--tree, $ \Sigma $, with finite edge stabilisers.
Moreover, each element of $ {\script H} $ is elliptic with respect to
$ \Sigma $, contradicting indecomposability.
\par
We thus have $ T = \bigcup {\script B} $.
Suppose, for contradiction, that there is an edge of $ T $ which lies
in precisely one element of $ {\script B} $.
We may as well suppose that this is true of all edges of $ T $.
(For if not, let $ F $ be the union of all edges of $ T $ which lie
in at least two elements of $ {\script B} $.
Collapsing each component of $ F $ to a point, we obtain a new
$ G $--tree.
We replace $ {\script B} $ by the set of axis of those elements of
$ {\script H} $ which remain hyperbolic.
Thus each element of the new arc system is the result of collapsing
an element of the old arc system along a collection of disjoint
compact subarcs.)
\par
We now construct a bipartite graph, $ \Sigma $, with vertex set
an abstract disjoint union of $ V(T) $ and $ {\script B} $, by deeming
$ x \in V(T) $ and $ \beta \in {\script B} $ to be adjacent in $ \Sigma $ if
$ x \in \beta $ in $ T $.
Now, it's easily verified that $ \Sigma $ is a simplicial tree, and that
the stabiliser of each pair $ (x,{\script B}) $ is finite.
In other words, $ \Sigma $ is a $ G $--tree with finite edge stabilisers.
Finally, we note that each element of $ {\script H} $ is elliptic in
$ \Sigma $.
This again contradicts the indecomposability of $ {\script H} $.
\endprf
\par
We now move on to considering splittings over two-ended subgroups.
Suppose that $ \Gamma $ is a group, and that $ \Sigma $ is a cofinite
$ \Gamma $--tree (with no terminal vertex) and with two-ended edge-stabilisers.
We can write $ V(\Sigma) $ as a disjoint union,
$ V(\Sigma) = V_1(\Sigma) \sqcup V_2(\Sigma) \sqcup V_\infty(\Sigma) $,
depending on whether the corresponding vertex stabiliser is one, two
or infinite-ended.
Note that $ V_2(\Sigma) $ is precisely the set of vertices of finite
degree.
\par
We remark that if there is a bound on the order of finite subgroups of
$ \Gamma $, and there are no infinitely divisible elements, then each
two-ended subgroup lies in a unique maximal two-ended subgroup.
In this case, we can refine our splitting so that for each vertex
$ v \in V_1(\Sigma) \cup V_\infty(\Sigma) $, the incident edge groups
are all maximal two-ended subgroups of $ \Gamma(v) $.
This is automatically true of the JSJ splitting of hyperbolic groups
(as described in [\Bob]), for example, though we shall have no need to
assume this in this section.
\par
It is fairly easy to see that the one-endedness or otherwise of $ \Gamma $
depends only on the infinite-ended vertex groups, $ \Gamma(v) $ for
$ v \in V_\infty(\Sigma) $.
In one direction, it easy to see that if one of these groups splits over a
finite group relative to incident edge groups, then we can refine our splitting
so that one of the new edge groups is finite.
Hence $ \Gamma $ is not one-ended.
In fact, we also have the converse.
Recall that a ``trivial vertex'' of a splitting is a vertex of
degree 1 such that the vertex group equals the adjacent edge group
(ie\ it corresponds to a terminal vertex of the corresponding tree).
\gap
\theorem{2.3}
{\sl
Suppose we represent a group, $ \Gamma $, as finite graph of
groups with two-ended vertex groups and no trivial vertices.
Then, $ \Gamma $ is one-ended if and only if none of the infinite-ended
vertex groups split intrinsically over a finite subgroup relative to
the incident edge groups.
}
\gap
\prf
Let $ \Sigma $ be the $ \Gamma $--tree corresponding to the splitting,
and write $ V(\Sigma) = V_1(\Sigma) \sqcup V_2(\Sigma) \sqcup
V_\infty(\Sigma) $ as above.
Given $ v \in V(\Sigma) $ let $ \Delta(v) \subseteq E(\Sigma) $
be the set of incident edges.
We are supposing that for each $ v \in V_\infty(\Sigma) $, the
set of incident edge stabilisers,
$ \{ \Gamma_\Sigma(e) \mid e \in \Delta(v) \} $, is indecomposable in
the group $ \Gamma_\Sigma(v) $.
This is therefore true for all $ v \in V(\Sigma) $.
We aim to show that $ \Gamma $ is one-ended.
\par
Suppose, for contradiction, that there exists a non-trivial minimal
$ G $--tree, $ T $, with finite edge stabilisers.
Let $ {\script B} $ be the arc system on $ T $ consisting of the axes
of those $ \Sigma $--edge stabilisers, $ \Gamma_\Sigma(e) $, which are
hyperbolic with respect to $ T $.
By Lemma 2.1, $ {\script B} $ is edge-finite.
\par
Suppose, first, that $ {\script B} = \emptyset $, ie\ each group
$ \Gamma_\Sigma(e) $ for $ e \in E(\Sigma) $ is elliptic in $ T $.
Suppose $ v \in V(\Sigma) $.
Since $ \{ \Gamma_\Sigma(e) \mid e \in \Delta(v) \} $ is indecomposable
in $ \Gamma_\Sigma(v) $, it follows that $ \Gamma_\Sigma(v) $ must
be elliptic in $ T $.
It therefore fixes a unique vertex of $ T $.
Suppose $ w \in V(\Sigma) $ is adjacent to $ v $.
Since $ \Gamma_\Sigma(v) \cap \Gamma_\Sigma(w) $ is infinite, it follows
that $ \Gamma_\Sigma(w) $ must also fix the same vertex of $ T $.
Continuing in this way, we conclude that this must be true of all
$ \Sigma $--vertex stabilisers.
We therefore arrive at the contradiction that $ \Gamma $ fixes a vertex of
$ T $.
\par
We deduce that $ {\script B} \ne \emptyset $.
Now, choose any $ \beta \in {\script B} $ and any edge $ \epsilon \in
E(\beta) $.
By construction, $ \beta $ is the axis of some edge stabiliser
$ \Gamma_\Sigma(e_0) $ for $ e_0 \in E(\Sigma) $.
Let $ v \in V(\Sigma) $ be an endpoint of $ e_0 $.
Now, $ \Gamma_\Sigma(e_0) \subseteq \Gamma_\Sigma(v) $, so $ \Gamma_\Sigma(v) $
is not elliptic in $ T $.
It follows that $ v \notin V_1(\Sigma) $.
If $ v \in V_2(\Sigma) $, then $ \beta $ is the axis in $ T $ of
$ \Gamma_\Sigma(v) $, and hence of any edge $ e_1 \in E(\Sigma) $
adjacent to $ e_0 $.
In particular, $ \epsilon $ lies in the axis of $ \Gamma_\Sigma(e_1) $.
If $ v \in V_\infty(\Sigma) $, let $ T(v) $ be the unique minimal
$ \Gamma_\Sigma(v) $--invariant subtree of $ T $.
Let $ {\script B}(v) $ be the set of axis of hyperbolic elements
of $ \{ \Gamma_\Sigma(e) \mid e \in \Delta(v) \} $.
Thus, $ {\script B}(v) \subseteq {\script B} $ is an arc system on
$ T(v) $, and $ \beta \in {\script B}(v) $.
By Lemma 2.2, there is some $ \beta' \in {\script B}(v) \bksl \{ \beta \} $
with $ \epsilon \in E(\beta') $.
Now, $ \beta' $ is the axis of $ \Gamma_\Sigma(e_1) $ for some edge
$ e_1 \in E(\Sigma) $ adjacent to $ e_0 $, as in the case where
$ v \in V_2(\Sigma) $.
Now, in the same way, we can find some edge $ e_2 $ incident on the
other endpoint of $ e_1 $, so that $ \Gamma_\Sigma(e_2) $ is
hyperbolic in $ T $ and contains $ \epsilon $ in its axis.
Continuing, we get an infinite sequence of edges, $ (e_n)_{n \in {\bbf N}} $,
which form a ray in $ \Sigma $, and which all have this property.
\par
Now, since $ {\script B} $ is edge-finite, we can pass to a subsequence
so that the axes of the groups $ \Gamma_\Sigma(e_n) $ are constant.
Since $ \Sigma $ is cofinite, we can find an edge $ e \in E(\Sigma) $
and an element $ \gamma \in \Gamma $ which is hyperbolic in $ \Sigma $,
and such that the axes of $ \Gamma_\Sigma(e) $ and
$ \Gamma_\Sigma(\gamma e) = \gamma \Gamma_\Sigma(e) \gamma^{-1} $
in $ T $ are equal to $ \alpha $, say.
In particular, $ \gamma \alpha = \alpha $.
Now, $ \Gamma_\Sigma(e) $ has finite index in the setwise stabiliser
of $ \alpha $, and so some power of $ \gamma $ lies in $ \Gamma_\Sigma(e) $,
contradicting the fact that $ \gamma $ is hyperbolic in $ \Sigma $.
\par
This finally contradicts the existence of the $ \Gamma $--tree $ T $.
\endprf
\par
We note that Theorem 2.3 gives a means of describing the indecomposibility
of a set of two-ended subgroups in terms of the ``doubled'' group, as follows.
\par
Suppose that $ G $ is a group, and that $ {\script H} $ is a union
of conjugacy classes of subgroups.
We form a graph of groups with two vertices as follows.
We take two copies of $ G $ as vertices, and connect them by a set
of edges, one for each conjugacy class of subgroup in $ {\script H} $.
We associate to each edge the corresponding group.
We refer to the fundamental group of this graph of groups as the \em{double} of
$ G $ in $ {\script H} $, and write it as $ D(G,{\script H}) $.
For example, if $ H $ is any subgroup of $ G $ and $ {\script H} $ is
its conjugacy class, then we just get the amalgamated free product,
$ D(G,{\script H}) \cong G \ast_H G $.
\par
 From Theorem 2.3, we deduce immediately:
\gap
\corol{2.4}
{\sl
Suppose that $ G $ is a group, and that $ {\script H} $ is a union
of finitely many conjugacy classes of two-ended subgroups.
Then, $ {\script H} $ is indecomposable in $ G $ if and only
if the double, $ D(G,{\script H}) $, is one-ended.
}
\endprf
\par
We note that Theorem 2.3 can be extended to allow for one-ended
edge groups.
The hypotheses remain unaltered.
We simply demand that no vertex group splits over a finite group
relative to the set of two-ended incident edge groups.
The argument remains essentially unchanged.
If, however, we allow for infinite-ended edge groups, then
Theorem 2.3 and Corollary 2.4 may fail.
\par
Consider, for example, a one-ended group, $ K $, with an infinite order
element $ a \in K $.
Let $ G $ be the free product $ K \ast {\bbf Z} $, and write $ b \in G $
for the generator of the $ {\bbf Z} $ factor.
Let $ H \le G $ be the subgroup generated by $ a $ and $ b $.
Thus, $ H $ is free of rank 2.
Now, the conjugacy class of $ H $ is indecomposable in $ G $.
(For suppose that $ T $ is a $ G $--tree with finite edge stabilisers
and with $ H $ elliptic.
Now, since $ K $ is one-ended, it is also elliptic.
Since $ K \cap H $ is infinite, and since $ K \cup H $ generates $ G $,
we arrive at the contradiction that $ G $ is elliptic.)
However, $ G \ast_H G $ is not one-ended.
In fact, $ G \ast_H G \cong (K \ast_{\langle a \rangle} K) \ast {\bbf Z} $.
We remark that by taking $ \langle a \rangle $ to be malnormal in $ K $
(for example taking $ K $ to be any torsion-free one-ended word hyperbolic
group, and taking $ a $ to be any infinite order element which is not a proper
power) we can arrange that $ H $ is malnormal in $ G $.
\gap\gap
\section{Indecomposable arc systems}
\gap\par
In this section, we look further at arc systems and give a combinatorial
characterisation of indecomposability.
First, we introduce some additional notation concerning trees.
\par
Suppose $ S \subseteq T $ is a subtree.
We write $ \pi_S\co T\cup \partial T \longrightarrow S \cup \partial S $
for the natural retraction.
Thus, $ \pi_S((T \cup \partial T) \bksl (S \cup \partial S))
\subseteq V(S) \subseteq S $.
If $ R \subseteq S $ is another subtree, then
$ \pi_R \circ \pi_S = \pi_R $.
Moreover, $ \pi_R|(S \cup \partial S) $ is defined intrinsically to
$ S $.
\par
If $ v \in V(S) $, then $ T \cap \pi_S^{-1}(v) $ is a subtree of
$ T $, which we denote by $ F(S,v) $.
Note that $ F(s,v) \cap S = \{ v \} $, and that
$ \partial F(S,v) = \partial T \cap \pi_S^{-1}(v) $.
Also, $ T = S \cup \bigcup_{v \in V(S)} F(S,v) $.
\par
We begin by describing generalisations of Whitehead graphs.
For the moment, we do not need to introduce group actions.
\par
Let $ T $ be a simplicial tree.
We write $ {\script S}(T) $ for the set of finite subtrees of $ T $.
We can think of $ {\script S}(T) $ as a directed set under inclusion.
Given $ S \in {\script S}(T) $, we define an equivalence relation,
$ \mathord{\approx}_S $, on $ \partial T $ by writing
$ x \approx_S y $ if $ \pi_S x = \pi_S y $.
In other words, $ x \approx_S y $ if and only if the arc connecting
$ x $ to $ y $ meets $ S $ in at most one point.
Clearly, if $ S \subseteq R \in {\script S}(T) $, then
$ \mathord{\approx}_R $ is finer than $ \mathord{\approx}_S $.
We therefore get a direct limit system of equivalence relations indexed
by $ {\script S}(T) $.
The direct limit (ie\ intersection) of these relations is just
the equality relation on $ \partial T $.
\par
Suppose now that $ {\script B} $ is an arc system on $ T $.
We have another equivalence relation, $ \mathord{\approx}_{\script B} $,
on $ \partial T $ defined as follows.
We write $ x \approx_{\script B} y $ if $ x=y $ or if there exists some
$ \beta \in {\script B} $ such that $ \partial \beta = \{ x,y \} $.
If the intersection of any two arcs of $ {\script B} $ is compact
(as in most of the cases in which we shall be interested) then
this is already an equivalence relation.
If not, we take $ \approx_{\script B} $ to be the transitive closure
of this relation.
\par
Given $ S \in {\script S}(T) $, let $ \mathord{\sim}_{S,{\script B}} $
be the transitive closure of the union of the relations
$ \mathord{\approx}_S $ and $ \mathord{\approx}_{\script B} $.
Thus, the relations $ \mathord{\sim}_{S,{\script B}} $ again form
a direct limit system indexed by $ {\script S}(T) $.
We write $ \mathord{\sim}_{\script B} $ for the direct limit.
\gap
\defne
We say that the arc system $ {\script B} $ is \em{indecomposable}
if there is just one equivalence class of $ \mathord{\sim}_{\script B} $
in $ \partial T $.
\gap
\par
We can give a more intuitive description of this construction which
ties in with Whitehead graphs as follows.
We fix our arc system $ {\script B} $.
If $ S \in {\script S}(T) $, we abbreviate
$ \mathord{\sim}_{S,{\script B}} $ to $ \mathord{\sim}_S $.
Note that, if $ Q \subseteq \partial T $ is a
$ \mathord{\sim}_S $--equivalence class, then
$ Q = \partial T \cap \pi_S^{-1} \pi_S Q $.
Let $ {\script W}(S) $ be the collection of all sets of the
form $ \pi_S Q $, as $ Q $ runs over the set,
$ \partial T/\mathord{\sim}_S $, of $ \mathord{\sim}_S $--classes.
Thus, $ {\script W}(S) $ gives a partition of the subset
$ \bigcup {\script W}(S) $ of $ V(S) $.
We refer to $ {\script W}(S) $ as a ``subpartition'' of
$ V(S) $ (ie\ a collection of disjoint subsets).
There is a natural bijection between $ {\script W}(S) $ and the
set $ \partial T/\mathord{\sim}_S $.
\par
Let us now suppose that $ \bigcup {\script B} $ is not contained in
any proper subtree of $ T $ (for example if $ {\script B} $ is
indecomposable).
Let $ {\script B}(S) \subseteq {\script B} $ be the set of arcs
which meet $ S $ in a non-trivial interval (ie\ non-empty and not a
point).
If $ \beta \in {\script B}(S) $, we write $ I(\beta) $ for the interval
$ \beta \cap S $, thought of abstractly, and write $ \fr I(\beta) $
for the set consisting of its two endpoints.
Let $ Z(S) $ be the disjoint union $ Z(S) =
\bigsqcup_{\beta \in {\script B}(S)} I(\beta) $, and let $ \fr Z(S) =
\bigsqcup_{\beta \in {\script B}(S)} \fr I(\beta) $.
There is a natural projection $ p\co Z(S)\longrightarrow S $ with
$ p(\fr Z(S)) \subseteq V(S) $.
Now let $ {\script G}(S) $ be the quotient space $ Z(S)/\mathord{\cong} $,
where $ \mathord{\cong} $ is the equivalence relation on $ Z(S) $
defined by $ x \cong y $ if and only if $ x=y $ or $ x,y \in \fr Z(S) $
and $ px = py $.
We see that $ {\script G}(S) $ is a 1--complex, with vertex set,
$ V({\script G}(S)) $, arising from $ \fr Z(S) $.
The map $ p $ induces a natural map from $ {\script G}(S) $ to
$ S $, also denoted by $ p $.
Now, $ p|V({\script G}(S)) $ is injective, and $ p(V({\script G}(S))) =
\bigcup {\script W}(S) $, where $ {\script W}(S) $ is the subpartition of
$ V(S) $ described earlier.
Moreover, an element of $ {\script W}(S) $ is precisely the vertex
set of connected component of $ {\script G}(S) $.
If $ {\script B} $ is edge-finite, then $ {\script G}(S) $ will
be a finite graph.
\par
To relate this to the theory of Whitehead graphs, the following
observation will be useful.
Recall that a graph is 2--vertex connected if it is connected and has
no cut vertex.
(We consider a graph consisting of a single edge to be 2--vertex connected.)
\gap
\lemma{3.1}
{\sl
Suppose that $ S_1, S_2 \in {\script S}(T) $ are such that
$ S_1 \cap S_2 $ consists of a single edge $ e \in E(S_1) \cap E(S_2) $.
If $ {\script G}(S_1) $ and $ {\script G}(S_2) $ are 2--vertex connected,
then so is $ {\script G}(S) $.
}
\gap
\prf
Let $ S = S_1 \cup S_2 \in {\script S}(T) $.
Let $ v_1, v_2 $ be the endpoints of $ e $ which are extreme in
$ S_1 $ and $ S_2 $ respectively.
Let $ V_1 = V(S_1) \bksl \{ v_1 \} $ and $ V_2 = V(S_2) \bksl \{ v_2 \} $.
Write $ W_i = p^{-1}(V_i) \subseteq V({\script G}(S)) $ so that
$ V({\script G}(S)) = W_1 \sqcup W_2 $.
Let $ {\script G}_i $ be the full subgraph spanned by $ W_i $.
Then $ {\script G}(S_i) $ is obtained by collapsing $ {\script G}_i $
to a single vertex.
The result therefore follows from the following observation, of which
we omit the proof.
\endprf
\gap
\lemma{3.2}
{\sl
Suppose that $ {\script G} $ is a connected graph and that
$ {\script G}_1 $ and $ {\script G}_2 $ are disjoint connected
subgraphs.
Write $ {\script G}'_i $ for the result of collapsing $ {\script G}_i $
to a single point in the graph $ {\script G} $.
If $ {\script G}'_1 $ and $ {\script G}'_2 $ are both 2--vertex connected,
then so is $ {\script G} $.}\break\hbox to 0pt{}\endprf
\par
Suppose $ v \in V(T) $.
Write $ S(v) $ for the subtree consisting of the union of all edges
incident on $ v $.
If $ T $ is locally finite, then $ S(v) \in {\script S}(T) $.
Applying Lemma 3.1 inductively we conclude:
\gap
\lemma{3.3}
{\sl
Suppose that $ {\script B} $ is an arc system on the locally finite
tree, $ T $, such that $ \bigcup {\script B} $ is not contained in any
proper subtree.
If $ {\script G}(S(v)) $ is 2--vertex connected for all $ v \in V(T) $,
then $ {\script B} $ is indecomposable.
}
\endprf
\par
The classical example of this, as discussed in the introduction, is
that of Whitehead graphs.
Suppose that $ G $ is a free group with free generators
$ a_1 ,\ldots, a_n $.
Let $ T $ be the Cayley graph of $ G $ with respect to these generators.
Thus, $ T $ is locally finite cofinite $ G $--tree.
\par
Let $ \{ \gamma_1 ,\ldots, \gamma_p \} $ be a finite set of non-trivial
elements of $ G $.
It's easy to see that the indecomposability of the set of cyclic
subgroups $ \{ \bet{\gamma_1} ,\ldots \bet{\gamma_p} \} $ (as defined in
Section 2) is equivalent to that of $ \{ H_1 ,\ldots, H_p \} $ where
$ H_k $ is the maximal cyclic subgroup containing $ \bet{\gamma_k} $.
For this reason, we don't loose any generality by taking the elements
$ \gamma_k $ to be indivisible, though this is not essential for what
are going to say.
\par
Now, let $ {\script B} $ be the arc system consisting of the
set of axes of all conjugates of the elements $ \gamma_i $.
Now, the graph $ {\script G}(S(v)) $ is independent of the choice
of vertex $ v \in V(T) $, so we may write it simply as $ {\script G} $.
We can construct $ {\script G} $ abstractly as the graph with
vertex set $ \{ a_1 ,\ldots, a_n, a_1^{-1} ,\ldots, a_n^{-1} \} $
where the number of edges connecting $ a_i^{\epsilon_i} $ to
$ a_j^{\epsilon_j} $ equals the total number of times the
subword $ a_i^{\epsilon_i} a_j^{-\epsilon_j} $ occurs in the (disjoint union of
the) reduced cyclic words representing elements $ \gamma_k $
(where $ \epsilon_i, \epsilon_j \in \{ -1,1 \} $).
Thus, the total number of edges in $ {\script G} $ equals the sum of
the cyclically reduced word lengths of the elements $ \gamma_k $.
The fact that we are taking reduced cyclic words tells us immediately
that there are no loops in $ {\script G} $.
We call $ {\script G} $ the \em{Whitehead graph}.
This agrees with the description in the introduction, except that we
are now allowing for multiple edges.
(To recover the description of the introduction, and that of the
original paper [\W], we can simply replace each multiple edge by a single edge.
This has no consequence for what we are going to say.)
\par
By Lemma 3.3, we see immediately that:
\gap
\propn{3.4}
{\sl
If $ {\script G} $ is 2--vertex connected, then $ {\script B} $ is
indecomposable.
}
\endprf
\par
We shall see later, in a more general context, that the
indecomposability of $ {\script B} $ is equivalent to the indecomposability
of the set of subgroups $ \{ \bet{\gamma_1} ,\ldots, \bet{\gamma_p} \} $.
\par
By a ``cut vertex'' of $ {\script G} $ we mean a vertex of $ {\script G} $
which separates the component in which it lies.
Now, if $ {\script G} $ contains a cut vertex, one can change the
generators (in an explicit algorithmic fashion) so as to reduce the
total length of $ {\script G} $ (allowing multiple edges) --- cf\ [\W].
Thus, after a linearly bounded number of steps, we arrive at a Whitehead
graph with no cut vertex.
(It follows that if we choose generators so as to minimise the sum of the
cyclically reduced word lengths of the $ \gamma_k $, then the Whitehead graph
will have this property.)
In this case, the Whitehead graph is either disconnected or 2--vertex connected.
In the former case, $ {\script B} $ is clearly not indecomposable, whereas
in the latter case it is (by Proposition 3.4).
There is therefore a linear algorithm to decide indecomposability for a finite
set of elements in a free group.
\par
We remark that we can also recognise a free generating set by the
same process.
If $ p=n $, then $ \{ \gamma_1 ,\ldots, \gamma_n \} $ forms a free generating
set if and only if a minimal Whitehead graph (or any Whitehead graph without
cut vertices) is a disjoint union of $ n $ bigons.
(If the elements $ \gamma_i $ are all indivisible, then any component with
2 vertices must be a bigon.)
The algorithm arising out of this procedure was one of the main motivations of
the original paper [\W].
\par
We want to generalise some of this discussion of indecomposability to
the context of groups accessible over finite groups, as alluded to
in Section 2.
\par
For the moment, suppose that $ G $ is any group, and that
$ T $ and $ \Sigma $ are equivalent cofinite $ G $--trees with
finite edge stabilisers.
There are morphisms $ \phi\co T\longrightarrow \Sigma $ and
$ \psi\co \Sigma \longrightarrow T $.
These morphisms are quasiisometries, and hence induce a canonical
bijection between $ \partial T $ and $ \partial \Sigma $.
In this case, it is appropriate to deal with formal arc systems,
ie\ ($ G $--invariant) sets of unordered pairs of elements of
$ \partial T \equiv \partial \Sigma $.
Such a formal arc system determines an arc system, $ {\script B} $, on
$ T $ and one, $ {\script A} $, on $ \Sigma $.
There is a bijection between $ {\script B} $ and $ {\script A} $ such
that corresponding arcs have the same ideal endpoints.
Thus, if $ \beta \in {\script B} $, then $ \phi(\beta) $ is a subtree
of $ \Sigma $, with $ \partial \phi(\beta) \equiv \partial \beta $.
We see that the corresponding arc, $ \alpha \in {\script A} $ is the unique
biinfinite arc contained in $ \phi(\beta) $.
Note that we get relations $ \sim_{\script B} $ and $ \sim_{\script A} $
on $ \partial T \equiv \partial \Sigma $, from the direct limit construction
described earlier.
Our first objective will be to check that these are equal.
It follows that the indecomposability of $ {\script A} $ and
$ {\script B} $ are equivalent (Lemma 3.5).
We thus get a well-defined notion of indecomposability of formal
arc systems for such trees.
\par
Suppose that $ S \in {\script S}(T) $.
For clarity, we write $ \mathord{\approx}_{S,T} $ for the relation on
$ \partial T $ abbreviated to $ \mathord{\approx}_S $ in the previous
discussion (ie\ $ x \approx_{S,T} y $ if $ \pi_S x = \pi_S y $).
We thus have a direct limit system
$ (\mathord{\approx}_{S,T})_{S \in {\script S}(T)} $.
We similarly get another direct limit system
$ (\mathord{\approx}_{R,\Sigma})_{R \in {\script S}(\Sigma)} $.
We claim that these are cofinal.
In other words, for each $ S \in {\script S}(T) $, there is some
$ R \in {\script S}(\Sigma) $ such that the relation
$ \mathord{\approx}_{R,\Sigma} $ is finer than
$ \mathord{\approx}_{S,T} $, and conversely, swapping the roles of
$ T $ and $ \Sigma $.
\par
To see this, let $ \phi \co T\longrightarrow \Sigma $ be a morphism,
and let $ T' $ be an equivariant subdivision of $ T $ such that
$ \phi\co T'\longrightarrow \Sigma $ is a folding.
Suppose $ R \in {\script S}(\Sigma) $.
Applying Lemma 1.2, there is finite subtree, $ S $, of $ T $ which contains
every edge of $ T' $ that gets mapped homeomorphically to one of the edges
of $ R $.
Suppose that $ x,y \in \partial T \equiv \partial \Sigma $, and let $ \alpha $
and $ \beta $ be the arcs in $ T $ and $ \Sigma $ respectively, connecting
$ x $ to $ y $.
Thus $ \beta \subseteq \phi \alpha $.
Suppose that $ x \approx_{S,T} y $.
In other words, $ \alpha \cap S $ is either empty or consists of a single
vertex.
We claim that the same is true of $ \beta \cap R $.
For any edge of $ \beta \cap R $ is the image under $ \phi $ of
some edge $ \epsilon $ of $ \alpha $ in $ T' $.
By construction, $ \epsilon $ is also an edge of $ S $ in $ T' $, giving
a contradiction.
This shows that $ x \approx_{R,\Sigma} y $ as claimed.
Swapping the roles of $ T $ and $ \Sigma $, we deduce the cofinality
of the direct limit systems as claimed.
\par
Now, suppose that $ {\script B} $ and $ {\script A} $ are arc
systems on $ T $ and $ \Sigma $ respectively, giving rise to
the same formal arc system.
We get identical relations $ \mathord{\approx}_{\script B} =
\mathord{\approx}_{\script A} $ on $ \partial T = \partial \Sigma $,
as defined earlier.
Now, it follows that the direct limit systems
$ (\mathord{\sim}_{S,{\script B}})_{S \in {\script S}(T)} $ and
$ (\mathord{\sim}_{R,{\script A}})_{R \in {\script S}(\Sigma)} $ are
cofinal, and so give rise to the same direct limit, namely
$ \mathord{\sim}_{\script B} = \mathord{\sim}_{\script A} $, as
claimed earlier.
\par
In particular, we see that $ {\script B} $ is indecomposable if
and only if $ {\script A} $ is.
In summary, reintroducing the group action, we have shown:
\gap
\lemma{3.5}
{\sl
Suppose that $ T $ and $ \Sigma $ are equivalent cofinite
$ G $--trees with finite edge stabilisers.
Suppose that $ {\script B} $ and $ {\script A} $
are arc systems on $ T $ and $ \Sigma $ respectively, corresponding
to the same formal arc system on $ \partial T \equiv \partial \Sigma $.
Then, $ {\script B} $ is indecomposable if and only if
$ {\script A} $ is indecomposable.
}
\endprf
\par
Suppose, now, that $ G $ is accessible over finite groups.
As discussed in Section 1, we can associate to $ G $ a set
$ \partial_\infty G $, which we can identify with the boundary
of any cofinite $ G $--tree with finite edge stabilisers and
finite and one-ended vertex stabilisers.
We refer to such trees as \em{complete} $ G $--trees.
Any two complete $ G $--trees are equivalent, so by Lemma 3.5, it makes sense
to speak about a formal arc system on $ \partial_\infty G $ as being
indecomposable.
\par
Suppose, now that $ H \le G $ is a two-ended subgroup.
We say that $ H $ is \em{elliptic} if it lies inside some one-ended
subgroup of $ G $.
Thus $ H $ is elliptic if and only it is elliptic with respect to
some (and hence any) complete $ G $--tree.
Otherwise, we say that $ H $ is \em{hyperbolic}.
In this case, there is a unique $ H $--invariant unordered pair of
points in $ \partial_\infty G $ which we denote by $ \Lambda H $.
Thus, $ \Lambda H $ is the pair of endpoints of the axis of $ H $
in any complete $ G $--tree.
We refer to $ \Lambda H $ as the \em{limit set} of $ H $.
We note that if $ H' $ is another hyperbolic two-ended subgroup, and
$ \Lambda H \cap \Lambda H' \ne \emptyset $, then $ H $ and $ H' $
are commensurable, and hence lie in the same maximal two-ended subgroup.
\par
Let $ {\script H} $ be a finite union of conjugacy classes of hyperbolic
two-ended subgroups of $ G $.
Recall that $ {\script H} $ is ``indecomposable'' if we cannot
write $ G $ as a non-trivial amalgamated free product or HNN--extension over
a finite group with each element of $ H $ conjugate into a vertex group.
It is easy to see that this property depends only on the
commensurability classes of the elements of $ {\script H} $, so we
may, if we wish, take all the elements of $ H $ to be maximal two-ended
subgroups, in which case their limit sets are all disjoint.
Note that we get a formal arc system,
$ \{ \Lambda H \mid H \in {\script H} \} $, on $ \partial_\infty G $.
We claim:
\gap
\propn{3.6}
{\sl
If the formal arc system $ \{ \Lambda H \mid H \in {\script H} \} $ is
indecomposable, then $ {\script H} $ is indecomposable.
}
\gap
\prf
Suppose not.
Then there is a non-trivial cofinite $ G $--tree, $ T $, with finite
edge stabilisers and with each element of $ {\script H} $ elliptic
with respect to $ T $.
Now, as discussed in Section 1, we can refine the splitting $ T/G $ to
a complete splitting, giving us a complete $ G $--tree, $ \Sigma $.
We can recover $ T $ by collapsing $ T $ along a disjoint union of
subtrees.
Each element of $ H $ fixes setwise one of these subtrees.
\par
Now, let $ {\script B} $ be the arc system on $ \Sigma $ given by the
formal arc system, in other words, the set of axes of elements of
$ {\script H} $.
Thus each axis lies inside one of the collapsing subtrees.
In particular, $ \Sigma \ne \bigcup {\script B} $, and so
$ {\script B} $ is decomposable.
\endprf
\par
We shall prove a converse to Proposition 3.6 in the case where
$ G $ is finitely generated.
For this we shall need a relative version of Stallings's theorem.
\par
Let $ G $ be a finitely generated group, and let $ X $ be a 
Cayley graph of $ X $ (or, indeed, any graph on which $ G $ acts with
finite vertex stabilisers and finite quotient).
Given a subset $ A \subseteq V(X) $ we write $ E_A \subseteq E(X) $
for the set of edges with precisely one endpoint in $ A $.
Thus, to say that $ X $ has ``more than one end'' means that we can
find an infinite subset, $ A \subseteq V(X) $ such that its complement
$ B = V(X) \bksl A $ is also infinite, and such that
$ E_A = E_B $ is finite.
Thus, Stallings's theorem [\St] tells us that in such a case, $ G $ splits
over a finite group.
\par
Suppose, now that $ H \le G $ is a two ended subgroup, and that
$ C \subseteq V(X) $ is an $ H $--orbit of vertices (or any $ H $--invariant
subset with $ C/H $ finite).
Now, for all but finitely many $ G $--images, $ gC $, of $ C $, we have either
$ gC \subseteq A $ or $ gC \subseteq B $.
For the remainder, we have three possibilities: either $ gC \cap A $ is
finite or
$ gC \cap B $ is finite, or else both of these subsets give us a neighbourhood
of an end of $ H $.
We shall not say more about the last case, since it is precisely the case we
wish to rule out.
Note that this classification does not depend on the choice of
$ H $--orbit, $ C $.
A specific relative version of Stallings's theorem says the following:
\gap
\lemma{3.7}
{\sl
Suppose $ G $ is a finitely group and $ {\script H} $ is a finite
union of conjugacy classes of two-ended subgroups.
Let $ X $ be a Cayley graph of $ G $.
Suppose we can find an infinite set, $ A \subseteq V(X) $, such
that $ E_A $ is finite and $ B = V(X) \bksl A $ is infinite.
Suppose that for any $ H \in {\script H} $ either $ A \cap C $ or
$ B \cap C $ is finite for some (hence every) $ H $--orbit of vertices,
$ C $.
Then, $ {\script H} $ is decomposable (ie\ $ G $ splits over a finite
group relative to $ {\script H} $).
}
\endprf
\par
In fact, a much stronger result follows immediately from the results of [\DiD].
It may be stated as follows.
Suppose $ G $ is any finitely generated group, and $ A \subseteq G $ is an
infinite subset, whose complement $ B = G \bksl A $ is also infinite.
Suppose that the symmetric difference of $ A $ and $ Ag $ is finite
for all $ g \in G $.
Suppose that $ H_1 ,\ldots, H_n $ are subgroups such
that for all $ g \in G $ and all $ i \in \{ 1 ,\ldots, n \} $ either
$ gH_i \cap A $ or $ gH_i \cap B $ is finite.
Then $ G $ splits over a finite group relative to $ \{ H_1 ,\ldots, H_n \} $.
(If fact, it's sufficient to rule out $ G $ being a non-finitely generated
countable torsion group.)
\par
Alternatively, one can deduce Lemma 3.7, as we have stated it, by applying
Stallings's theorem to the double, $ D(G,{\script H}) $, and
using Corollary 2.4.
We briefly sketch the argument.
We may construct a Cayley graph, $ Y $, for $ D(G,{\script H}) $ by
taking lots of copies of $ X $, and stringing them together in a
treelike fashion.
Let's focus on a particular copy of $ X $, which we take to be acted upon
by $ G $.
Now each adjacent copy of $ X $ corresponds to an element $ H \in {\script H} $,
and is connected ours by an $ H $--orbit of edges.
We refer to such edges as ``amalgamating edges''.
The amalgamating edges corresponding to $ H $ are attached to $ X $ by
an $ H $--orbit, $ C_H $, of vertices of $ X $.
By hypothesis, either $ C_H \cap A $ is finite, in which case, we
write $ E_H $ for the set of amalgamating edges which have an endpoint
in $ C_H \cap A $, or else, $ C_H \cap B $ is finite, in which case, we
write $ E_H $ for the set of amalgamating edges which have an endpoint
in $ C_H \cap B $.
Now, for all but finitely many $ H $, the set $ E_H $ is empty.
Thus, the set $ E_{\script H} = \bigcup_{H \in {\script H}} E_H $ is finite,
and so $ E_0 = E_A \cup E_{\script H} \subseteq E(Y) $ is finite.
Now, $ E_0 $ separates $ Y $ into two infinite components.
Thus, by Stallings's theorem, $ D(G,{\script H}) $ splits over a finite group,
and so by Corollary 2.4, $ {\script H} $ is decomposable.
With the details filled in, this gives another proof of Lemma 3.7.
\par
We are now ready to prove a converse to Proposition 3.6:
\gap
\propn{3.8}
{\sl
Suppose that $ G $ is a finitely generated accessible group.
Suppose that $ {\script H} $ is a finite union of conjugacy classes of
hyperbolic two-ended subgroups.
If $ {\script H} $ is indecomposable, then the formal arc system,
$ \{ \Lambda H \mid H \in {\script H} \} $, on $ \partial_\infty G $,
is indecomposable.
}
\gap
\prf
Let $ T $ be a complete $ G $--tree, and let $ {\script B} $ be
the corresponding arc system on $ T $, ie\ the set of axes of
elements of $ {\script H} $.
Suppose, for contradiction, that $ {\script B} $ is decomposable.
In other words, we can find $ S \in {\script S}(T) $ such that there
is more than one $ \mathord{\sim}_S $--class.
By taking projections of $ \mathord{\sim}_S $--classes as disussed in
Section 1, we can write $ V(S) $ as a disjoint union of non-empty
subsets, $ V(S) = W_1 \sqcup W_2 $ with the property that if
$ \beta \in {\script B} $, then $ \beta $ meets $ S $, if at all, in
compact interval (or point) with either both endpoints in $ W_1 $
or both endpoints in $ W_2 $.
Let $ F_i = \pi_S^{-1} W_i $.
Thus, $ T = S \cup F_1 \cup F_2 $, and each component of each $ F_i $
is a subtree meeting $ S $ in a single point.
\par
Now, let $ X $ be a Cayley graph of $ G $.
Let $ f\co V(X)\longrightarrow V(T) $ be any $ G $--equivariant map.
Let $ A_i = f^{-1} F_i \subseteq V(X) $.
Thus, $ V(X) = A_1 \sqcup A_2 $.
Moreover, it is easily seen that $ E_{A_1} = E_{A_2} $ is finite.
(For example, extend $ f $ equivariantly to a map $ f\co X\longrightarrow T $
so that each edge of $ X $ gets mapped to a compact interval of $ T $.
Only finitely many $ G $--orbits of such an interval can contain a given
edge of $ T $.
Now, the image of an edge of $ E_{A_1} $ connects a vertex of $ F_1 $ to
a vertex of $ F_2 $, and hence contains an edge of $ S $.
There are only finitely many such edges.)
\par
Finally, suppose that $ H \in {\script H} $.
Let $ \beta \in {\script B} $ be the axis of $ H $.
Without loss of generality, we can suppose that both ends of $ \beta $
are contained in $ F_1 $.
Now suppose that $ C $ is any $ H $--orbit of vertices of $ X $.
Then $ f(C) $ remains within a bounded distance of $ \beta $, from
which we see easily that $ f(C) \cap F_2 $ is finite.
Thus, $ C \cap A_2 $ is finite.
\par
We have verified the hypotheses of Lemma 3.7, and so $ {\script H} $ is
decomposable, contrary to our hypotheses.
\endprf
\par
Note that Propositions 3.6 and 3.8 apply, in particular, to any finitely
presented group, and even more specifically, to any hyperbolic group, $ G $.
In the latter case, $ \partial_\infty G $ can be identified as a subset
of the Gromov boundary, $ \partial G $, as discussed in Section 2.
If $ H \le G $ is a hyperbolic two-ended subgroup, then
$ \Lambda H \subseteq \partial G $ is the limit set of $ H $ by the
standard definition.
This ties in with the discussion of equivalence relations on
$ \partial G $ in the introduction, and will be elaborated on in
Section 5.
\gap\gap
\section{Quasiconvex splittings of hyperbolic groups}
\gap\par
For most of the rest of this paper, we shall be confining our attention
to hyperbolic groups.
We shall consider how some of the general constructions of Sections 1--3
relate to the topology of the boundary in this case.
Before we embark on this, we review some general facts about quasiconvex
splittings of hyperbolic groups (ie\ splittings over quasiconvex subgroups).
This elaborates on the account given in [\Bob].
\par
Throughout the rest of this paper, we shall use the notation
$ \fr A $ to denote the topological boundary (or ``frontier'') of
a subset, $ A $, of a larger topological space.
We reserve the symbol ``$ \partial $'' for ideal boundaries.
\par
Let $ \Gamma $ be any hyperbolic group.
Let $ X $ be any locally finite connected graph on which $ \Gamma $ acts
freely and cocompactly (for example a Cayley graph of $ \Gamma $).
We put a path metric, $ d $, on $ X $ by assigning a positive
length to each edge in a $ \Gamma $--invariant fashion.
Let $ \partial \Gamma \equiv \partial X $ be the boundary of $ \Gamma $.
We may put a metric on $ \partial \Gamma $ as described in [\GhH].
This has the property that given a basepoint, $ a \in V(X) $, there
are constants, $ A,B > 0 $ and $ \lambda \in (0,\infty) $ such that
if $ x,y \in \partial X $, then $ A \lambda^\delta \le \rho(x,y) \le
B \lambda^\delta $, where $ \delta $ is the distance from $ a $ to
some biinfinite geodesic connecting $ x $ to $ y $.
Although all the arguments of this paper can be expressed in purely
topological terms, it will be convenient to have recourse to this metric.
\par
Note that if $ G \le \Gamma $ is quasiconvex, then it is intrinsically
hyperbolic, and we may identify its boundary, $ \partial G $, with
its limit set $ \Lambda G \subseteq \partial \Gamma $.
Note that $ G $ acts properly discontinuously on $ \partial \Gamma \bksl
\Lambda G $.
The setwise stabiliser of $ \Lambda G $ in $ \Gamma $ is precisely
the commensurator, $ \Comm(G) $, of $ G $ in $ \Gamma $
(ie\ the set of all $ g \in \Gamma $ such that $ G \cap gGg^{-1} $
has finite index in $ G $).
In this case, $ G $ has finite index in $ \Comm(G) $.
In fact, $ \Comm(G) $ is the unique maximal subgroup of $ \Gamma $
which contains $ G $ as finite index subgroup.
We say that $ G $ is {\it full \/} if $ G = \Comm(G) $.
\par
We shall use the following notation.
If $ f\co Z\longrightarrow [0,\infty) $ is a function from some set
$ Z $ to the nonnegative reals, we write ``$ f(z) \rightarrow 0 $
for $ z \in Z $'' to mean
that $ \{ z \in Z \mid f(z) \ge \epsilon \} $ is finite for all
$ \epsilon > 0 $.
We similarly define ``$ f(z) \rightarrow \infty $ for $ z \in Z $''.
\gap
\lemma{4.1}
{\sl
If $ G \le \Gamma $ is quasiconvex and $ x \in \partial \Gamma $,
then $ \rho(gx, \Lambda G) \rightarrow 0 $ for $ g \in G $.
}
\gap
\prf
Since $ G $ acts properly discontinuously on $ \partial \Gamma \bksl
\Lambda G $, there can be no accumulation point of the $ G $--orbit
of $ x $ in this set.
\endprf
\par
The following is also standard:
\gap
\lemma{4.2}
{\sl
If $ G \le \Gamma $ is quasiconvex, then $ \diam(\Lambda H) \rightarrow 0 $
as $ H $ ranges over conjugates of $ G $.
}
\endprf
\par
We want to go on to consider splittings of $ \Gamma $.
For this, we shall want to introduce some further notation regarding trees.
\par
By a ``directed edge'' we mean an edge together with an orientation.
We write $ {\vec E}(T) $ for the set of directed edges.
We shall always use the convention that $ e \in E(T) $ represents the
undirected edge underlying the directed edge $ {\vec e} \in {\vec E}(T) $.
We write $ \rhead({\vec e}) $ and $ \rtail({\vec e}) $ respectively for
the head and tail of $ {\vec e} $.
We use $ -{\vec e} $ for the same edge oriented in the opposite
direction, ie\ $ \rhead(-{\vec e}) = \rtail({\vec e}) $ and
$ \rtail(-{\vec e}) = \rhead({\vec e}) $.
If $ {\vec e} \in {\vec E}(T) $ and $ v \in V(T) $, we say that
$ {\vec e} $ ``points towards'' $ v $ if $ \dist(v,\rtail({\vec e}))
= \dist(v,\rhead({\vec e}))+1 $.
\par
If $ v \in V(T) $, let $ \Delta(v) \subseteq E(T) $ be the set of
edges incident on $ v $, and let $ {\vec \Delta}(v) =
\{ {\vec e} \in {\vec E}(T) \mid \rhead({\vec e}) = v \} $.
Thus, the degree of $ v $ is $ \card(\Delta(v)) =
\card({\vec \Delta}(v)) $.
\par
Given $ {\vec e} \in {\vec E}(T) $, we write $ \Phi({\vec e}) =
\Phi_T({\vec e}) $ for the connected component of $ T $ minus the
interior of $ e $ which contains $ \rtail({\vec e}) $.
Thus, $ V(\Phi({\vec e})) $ is the set of vertices, $ v $, of
$ T $ such that $ {\vec e} $ points away from $ v $.
\par
Given $ v \in V(T) $, we shall write $ {\vec \Omega}(v) \subseteq
{\vec E}(T) $ for the set of directed edges which point towards $ v $.
Thus, for each edge $ e \in E(T) $, precisely one of the pair
$ \{ {\vec e}, -{\vec e} \} $ lies in $ {\vec \Omega}(v) $.
Note that $ {\vec e} \in {\vec \Omega}(v) $ if and only if
$ v \notin \Phi({\vec e}) $.
Clearly $ {\vec \Delta}(v) \subseteq {\vec \Omega}(v) $.
\par
We now return to our hyperbolic group, $ \Gamma $.
Suppose that $ \Gamma $ acts without edge inversions on a simplicial
tree, $ \Sigma $, with $ \Sigma/\Gamma $ finite.
We suppose that this action is minimal.
Given $ v \in V(\Sigma) $ and $ e \in E(\Sigma) $, write $ \Gamma(v) $
and $ \Gamma(e) $ respectively for the corresponding vertex and
edge stabilisers.
Note that $ \Gamma(v) $ is finite if and only if $ v $ has finite degree
in $ \Sigma $ and finite incident edge stabilisers.
If $ v,w \in V(\Sigma) $ are the endpoints of an edge $ e \in E(\Sigma) $,
then $ \Gamma(e) = \Gamma(v) \cap \Gamma(w) $.
\par
As in [\Bob], we may construct a $ \Gamma $--equivariant map
$ \phi\co X\longrightarrow \Sigma $ such that each edge of $ X $
either gets collapsed onto a vertex of $ \Sigma $ or mapped
homeomorphically onto a closed arc in $ \Sigma $.
(Note that, after subdividing $ X $ if necessary, we can assume that,
in the latter case, this closed arc is an edge of $ \Sigma $.)
Since the action of $ \Gamma $ is minimal, $ \phi $ is surjective.
\par
A proof of the following result can be found in [\Bob], though it
appears to be ``folklore''.
\gap
\propn{4.3}
{\sl
If $ \Gamma(e) $ is quasiconvex for each $ e \in E(\Sigma) $, then
$ \Gamma(v) $ is quasiconvex for each $ v \in V(\Sigma) $.
}
\endprf
\par
We refer to such a splitting as a {\it quasiconvex splitting\/}.
\par
We note that if a vertex group, $ \Gamma(v) $, of a quasiconvex
splitting has the property that all incident edge groups are of
infinite index in $ \Gamma(v) $, then $ \Gamma(v) $ must be full in
the sense described above.
In other words, $ \Gamma(v) $ is the setwise stabiliser of
$ \Lambda \Gamma(v) $.
This will be the case in most situations of interest (in particular
where all edge groups are finite or two-ended, but $ \Gamma(v) $ is not).
\par
Note that, if $ v,w \in V(\Sigma) $, then $ \Gamma(v) \cap \Gamma(w) $
is quasiconvex (since the intersection of any two quasiconvex subgroups
is quasiconvex [\Sh]).
We see that $ \Lambda \Gamma(v) \cap \Lambda \Gamma(w) =
\Lambda(\Gamma(v) \cap \Gamma(w)) $.
In particular, if $ v,w $ are the endpoints of an edge $ e \in E(\Sigma) $,
then $ \Lambda \Gamma(v) \cap \Lambda \Gamma(w) = \Lambda \Gamma(e) $.
\par
As described in [\Bob], there is a natural $ \Gamma $--invariant
partition of $ \partial \Gamma $ as $ \partial \Gamma =
\partial_0 \Gamma \sqcup \partial_\infty \Gamma $, where
$ \partial_0 \Gamma = \bigcup_{v \in V(\Sigma)} \Lambda \Gamma(v) $, and
$ \partial_\infty \Gamma $ is naturally identified with
$ \partial \Sigma $.
Note that $ \partial_\infty \Gamma $ is dense in $ \partial \Gamma $,
provided that $ \Sigma $ is non-trivial.
(In the case where the edge stabilisers are all finite, this
agrees with the notion introduced for accessible groups in Section 2.)
\par
Given $ {\vec e} \in  {\vec E}(\Sigma) $, we write
$$ \Psi({\vec e}) = \partial \Phi({\vec e}) \cup
\bigcup_{v \in V(\Phi({\vec e}))} \Lambda \Gamma(v) .$$
It's not hard to see that $ \Psi({\vec e}) $ is a closed
$ \Gamma(e) $--invariant subset of $ \partial \Gamma $. 
Moreover, $ \Psi({\vec e}) \cup \Psi(-{\vec e}) = \partial \Gamma $
and $ \Psi({\vec e}) \cap \Psi(-{\vec e}) =
\fr \Psi({\vec e}) = \Lambda \Gamma(e) $.
\par
Now, $ V(\Sigma) = \{ v \} \sqcup
\bigsqcup_{{\vec e} \in {\vec \Delta}(v)} V(\Phi({\vec e})) $ and
$ \partial \Sigma = \bigsqcup_{{\vec e} \in {\vec \Delta}(v)}
\partial \Phi({\vec e}) $.
It follows that:
\gap
\lemma{4.4}
{\sl
$ \partial \Gamma = \Lambda \Gamma(v) \cup
\bigcup_{{\vec e} \in {\vec \Delta}(v)} \Psi({\vec e}) $.
}
\endprf
\par
Moreover, for each $ {\vec e} \in {\vec \Delta}(v) $, we have
$ \Lambda \Gamma(v) \cap \Psi({\vec e}) = \Lambda \Gamma(e) $.
\par
The above assertions become more transparent, given the following
alternative description of $ \Psi({\vec e}) $.
\par
Let $ m(e) $ be the midpoint of the edge $ e $, and let $ I({\vec e}) $
be the closed interval in $ \Sigma $ consisting of the segment
of $ e $ lying between $ m(e) $ and $ \rtail({\vec e}) $.
Let $ Q(e) = \phi^{-1}(m(e)) \subseteq X $ and
$ R({\vec e}) = \phi^{-1}(\Phi({\vec e}) \cup I({\vec e})) \subseteq
X $, where $ \phi\co X\longrightarrow \Sigma $ is the map described above.
Note that $ Q(e) = \fr R({\vec e}) = R({\vec e}) \cap R(-{\vec e}) $.
By the arguments given in [\Bob], we see easily that $ Q(e) $ and
$ R({\vec e}) $ are quasiconvex subsets of $ X $.
Moreover, $ \Psi({\vec e}) = \partial R({\vec e}) $.
\par
Note that the collection $ \{ Q(e) \mid e \in E(\Sigma) \} $ is locally
finite in $ X $.
It follows that, for any fixed $ a \in X $, we have
$ d(a,Q(e)) \rightarrow \infty $ for $ e \in E(\Sigma) $.
\par
Now, fix some vertex, $ v \in V(\Sigma) $.
Recall that $ {\vec \Omega}(v) $ is defined to be the set
of all directed edges pointing towards $ v $.
Choose any $ b \in \phi^{-1}(v) \subseteq X $.
Now, if $ {\vec e} \in {\vec \Omega}(v) $, we have $ v \notin
\Phi({\vec e}) \cup I({\vec e}) $, and so $ b \notin R({\vec e}) $.
Since $ Q(e) = \fr R({\vec e}) $, we have $ d(b,R({\vec e})) =
d(b,Q(e)) $.
It follows that $ d(b,R({\vec e})) \rightarrow \infty $ for
$ {\vec e} \in {\vec \Omega}(v) $.
In fact, we see that $ d(a,R({\vec e})) \rightarrow \infty $ given
any fixed basepoint, $ a \in X $.
Now, there are only finitely many $ \Gamma $--orbits of directed
edges, and so the sets $ R({\vec e}) $ are uniformly quasiconvex.
 From the definition of the metric $ \rho $ on $ \partial \Gamma $,
it follows easily that $ \diam(\Psi({\vec e})) \rightarrow 0 $,
where $ \diam $ denotes diameter with respect to $ \rho $.
In summary, we have shown:
\gap
\lemma{4.5}
{\sl
For any $ v \in V(\Sigma) $, $ \diam(\Psi({\vec e})) \rightarrow 0 $
for $ {\vec e} \in {\vec \Omega}(v) $.
}
\endprf
\par
We now add the hypothesis that $ \Gamma(e) $ is infinite for all
$ e \in E(\Sigma) $.
\par
Suppose $ v \in V(\Sigma) $ and suppose $ K $ is any closed subset
of $ \Lambda \Gamma(v) $.
Let $ {\vec \Delta}_K(v) = \{ {\vec e} \in {\vec \Delta}(v) \mid
\Lambda \Gamma(e) \subseteq K \} $, and let
$ \Upsilon(v,K) = K \cup \bigcup_{{\vec e} \in {\vec \Delta}_K(v)}
\Psi({\vec e}) \subseteq \partial \Gamma $.
\gap
\lemma{4.6}
{\sl
The set $ \Upsilon(v,K) $ is closed in $ \partial \Gamma $.
}
\gap
\prf
Suppose $ x \notin \Upsilon(v,K) $.
In particular, $ x \notin K $, so $ \epsilon = \rho(x,K) >0 $.
Now, if $ {\vec e} \in {\vec \Delta}_K(v) $ and
$ \rho(x,\Psi({\vec e})) < \epsilon/2 $, then
$ \diam(\Psi({\vec e})) > \epsilon/2 $
(since $ K \cap \Psi({\vec e}) \supseteq \Lambda \Gamma(e) $, which,
by the hypothesis on edge stabilisers, is non-empty).
By Lemma 4.5, this occurs for only finitely many such $ {\vec e} $.
Since each $ \Psi({\vec e}) $ is closed, it follows that
$ \rho(x,\Upsilon(v,K)) $ is attained, and hence positive.
In other words, $ x \notin \Upsilon(v,K) $ implies
$ \rho(x,\Upsilon(v,K)) > 0 $.
This shows that $ \Upsilon(v,K) $ is closed.
\endprf
\gap\gap
\section{Quotients}
\gap\par
In this section, we aim to consider quotients of boundaries of
hyperbolic groups, and to relate this to indecomposability,
thereby generalising some of the results of [\O].
\par
First, we recall a few elementary facts from point-set topology [\Ke,\HoY].
Let $ M $ be a hausdorff topological space.
A subset of $ M $ is {\it clopen\/} if it is both open and closed.
We may define an equivalence relation on $ M $ by deeming two points
to be related if every clopen set containing one must also contain the
other.
The equivalence classes are called {\it quasicomponents\/}.
A {\it component\/} of $ M $ is a maximal connected subset.
Components and quasicomponents are always closed.
Every component is contained in a quasicomponent, but not conversely in
general.
However, if $ M $ is compact, these notions coincide.
Thus, if $ K $ and $ K' $ are distinct components of a compact hausdorff
space, $ M $, then there is a clopen subset of $ M $ containing $ K $, but
not meeting $ K' $.
\par
Suppose that $ M $ is a compact hausdorff space, and that
$ \mathord{\approx} $ is an equivalence relation on $ M $.
If the relation $ \mathord{\approx} $ is closed (as a subset of
$ M \times M $), then the quotient space, $ M/\mathord{\approx} $
is hausdorff.
\par
The compact spaces of interest to us here will be the boundaries of
hyperbolic groups.
Suppose that $ G $ is a hyperbolic group, and that $ \partial G $ is
its boundary.
Now, any two ended subgroup, $ H $, of $ G $ is necessarily quasiconvex,
so its limit set, $ \Lambda H \subseteq \partial G $, consists of pair
of points.
If $ H' $ is another two-ended subgroup, and $ \Lambda H \cap \Lambda H' 
\ne \emptyset $, then $ H $ and $ H' $ are commensurable, and so lie
in a common maximal two-ended subgroup.
In particular, $ \Lambda H = \Lambda H' $ (cf\ the discussion of accessible
groups in Section 3).
\par
Suppose that $ {\script H} $ is a union of finitely many conjugacy classes
of two-ended subgroups of $ G $.
Let $ {\mathord \approx}_{\script H} $ be the equivalence relation defined
on $ \partial G $ defined by $ x \approx_{\script H} y $ if and only if
either $ x=y $ or there exists $ H \in {\script H} $ such that
$ \Lambda H = \{ x,y \} $.
Now, it's a simple consequence of Lemma 4.2 that the relation
$ {\mathord \approx}_{\script H} $ is closed.
We write $ M(G,{\script H}) $ for the quotient space
$ \partial G/{\mathord \approx}_{\script H} $.
Thus:
\gap
\lemma{5.1}
{\sl
$ M(G,{\script H}) $ is compact hausdorff.
}
\endprf
\par
We aim to describe when $ M(G,{\script H}) $ is connected.
Clearly, if $ G $ is one-ended so that $ \partial G $ is connected,
this is necessarily the case.
We can thus restrict attention to the case when $ G $ is infinite-ended.
\par
Let $ T $ be a complete $ G $--tree.
As in Section 3, we can define $ \partial_\infty G $ as $ \partial T $.
This also agrees with the notation introduced in Section 4, thinking
of $ T $ as a quasiconvex splitting of $ G $.
In particular, we can identify $ \partial_\infty G $ as a subset of
$ \partial G $.
This set $ \partial_0 G = \partial G \bksl \partial_\infty G $
is a disjoint union of the boundaries of the infinite vertex stabilisers
of $ T $, ie\ the maximal one-ended subgroups.
In other words, the components of $ \partial_0 G $ are precisely the
boundaries of the maximal one-ended subgroups of $ G $.
\par
Let $ {\script H} $ be a set of two-ended subgroups as above.
The subset, $ {\script H}_0 $, of $ {\script H} $ consisting
of those subgroups in $ {\script H} $ which are hyperbolic (ie\ with
both limit points in $ \partial_\infty G $), defines a formal arc system
on $ \partial_\infty G $.
We aim to show that $ M(G,{\script H}) $ is connected if and only if
this arc system is indecomposable.
This, in turn, we know to be equivalent to asserting that $ {\script H}_0 $
is irreducible.
\par
In fact, it's easy to see that the elliptic elements of $ {\script H} $ have
no bearing on the connectivity or otherwise of $ M(G,{\script H}) $.
For this reason, we may as well suppose, for simplicity, that
$ {\script H} $ consists entirely of hyperbolic two-ended subgroups.
We therefore aim to show:
\gap
\theorem{5.2}
{\sl
Let $ G $ be an infinite-ended hyperbolic group, and let
$ {\script H} $ be a union of finitely many conjugacy classes
of hyperbolic two-ended subgroups.
Then, the quotient space $ M(G,{\script H}) $ is connected if and
only if $ {\script H} $ is indecomposable.
}
\gap
\ppar
First, we set about proving the ``only if'' bit.
Let $ T $ be a complete $ G $--tree.
Thus, $ \partial_\infty G $ is identified with $ \partial T $, and
$ {\script H} $ determines an arc system, $ {\script B} $, on $ T $.
We know (Propositions 3.6 and 3.8) that the indecomposability of
$ {\script H} $ is equivalent to the indecomposability of $ {\script B} $.
\par
We shall say that a subgraph, $ F $, of $ T $ is \em{finitely separated}
if there are only finitely many edges of $ T $ with precisely one
endpoint in $ F $.
Now, it's not hard to see that $ F $ is finitely separated if and only if
it's a finite union of finite intersections of subtrees of the form
$ \Phi({\vec e}) $ for $ {\vec e} \in {\vec E}(T) $ (recalling the
notation of Section 4).
\par
Now, given a subgraph, $ F \subseteq T $, we write
$$ A(F) = \partial F \cup \bigcup_{v \in V(F)} \Lambda G(v) $$
(so that $ A(T) = \partial G $).
If $ F $ is finitely separated, then $ A(F) $ is a finite union
of finite intersections of sets of the form $ \Psi({\vec e}) $, which
are each closed by the remarks of Section 4.
We conclude:
\gap
\lemma{5.3}
{\sl
If $ F \subseteq T $ is a finitely separated subgraph, then
$ A(F) $ is closed in $ \partial G $.
}
\endprf
\par
We can now prove:
\gap
\lemma{5.4}
{\sl
If $ M(G,{\script H}) $ is connected, then the arc system $ {\script B} $ is
indecomposable.
}
\gap
\prf
Suppose, to the contrary, that $ {\script B} $ is decomposable.
Then, exactly as in the proof of Proposition 3.8, we can find two disjoint
finitely separated subgraphs, $ F_1 $ and $ F_2 $ of $ T $
with $ V(T) = V(F_1) \sqcup V(F_2) $ and $ \partial T = \partial F_1 \sqcup
\partial F_2 $, and such that for each $ \beta \in {\script B} $, either
$ \partial \beta \subseteq \partial F_1 $ or
$ \partial \beta \subseteq \partial F_2 $.
We see that $ \partial G = A(F_1) \sqcup A(F_2) $.
\par
Let $ q\co\partial G \longrightarrow \partial G/{\mathord \approx}_{\script H} =
M(G,{\script H}) $ be the quotient map.
Now, from the construction, we see that if $ x \approx_{\script H} y $
then either $ x,y \in \partial F_1 \subseteq A(F_1) $ or
$ x,y \in \partial F_2 \subseteq A(F_2) $.
We therefore get that $ M(G,{\script H}) = q(A(F_1)) \sqcup q(A(F_2)) $.
But applying Lemma 5.3, the sets $ q(A(F_i)) $ are both closed in
$ M(G,{\script H}) $, contrary to the assumption that $ M(G,{\script H}) $ is
connected.
\endprf
\gap
\lemma{5.5}
{\sl
If $ {\script H} $ is indecomposable, then $ M(G,{\script H}) $ is
connected.
}
\gap
\prf
Suppose, for contradiction, that we can write $ M(G,{\script H}) $ as
the disjoint union of two non-empty closed sets, $ K_1 $ and $ K_2 $.
Let $ L_i \subseteq \partial G $ be the preimage of $ K_i $ under
the quotient map $ \partial G \longrightarrow M(G,{\script H}) $.
Thus, $ \partial G = L_1 \sqcup L_2 $.
Let $ X $ be a Cayley graph of $ G $.
Now, we can give $ X \cup \partial G $ a natural $ G $--invariant topology
as a compact metrisable space.
Since $ X \cup \partial G $ is normal, we can find disjoint open
subsets, $ U_i \subseteq X \cup \partial G $ with $ L_i \subseteq U_i $.
Now, $ (X \cup \partial G) \bksl (U_1 \cup U_2) \subseteq X $ is compact,
and so lies inside a finite subgraph, $ Y $, of $ X $.
Let $ A = U_1 \cap V(X) $ and let $ B = V(X) \bksl A $.
We need to verify that $ A $ satisfies the hypotheses of Lemma 3.7.
\par
Note that $ A \cup L_1 $ and $ B \cap L_2 $ are both closed in $ X \cup
\partial G $.
We see that $ A $ and $ B $ are both infinite.
Recall that $ E_A = E_B $ is the set of edges of $ X $ which have one
endpoint in $ A $ and the other in $ B $.
Now, $ E_A \subseteq E(Y) $, and so $ E_A $ is finite.
\par
Finally, suppose that $ H \in {\script H} $ and that $ C \subseteq V(X) $
is an $ H $--orbit of vertices of $ X $.
Now, $ C \cup \partial H $ is closed in $ X \cup \partial G $.
Without loss of generality we can suppose that $ \Lambda H \subseteq L_1 $.
Since $ B \cup L_2 \subseteq X \cup \partial G $ is closed, we see
that $ C \cap B $ is finite.
\par
We have verified the hypotheses of Lemma 3.7, and so we arrive at the
contradiction that $ {\script H} $ is decomposable.
\endprf
\par
This concludes the proof of Theorem 5.2.
\gap\gap
\section{Splittings of hyperbolic groups over finite and two-ended subgroups}
\gap\par
Suppose that a hyperbolic group splits over a collection of
two-ended subgroups.
We may in turn try to split each of the vertex groups over finite groups,
thus giving us a two-step series of splittings.
We want to study how the combinatorics of such splittings are reflected
in the topology of the boundary.
The combinatorics can be described in terms of the trees associated to
each step of the splitting, together with arc systems on the trees
of the second step which arise from the incident edge groups of
the first step.
\par
Suppose that $ \Gamma $ is a hyperbolic group, and that
$ \Sigma $ is a cofinite $ \Gamma $--tree with two-ended edge stabilisers.
Note that this is necessarily a quasiconvex splitting (as described in
Section 4), since a two-ended subgroup of a hyperbolic group is necessarily
quasiconvex (see, for example, [\GhH]).
We shall fix some vertex, $ \omega \in V(\Sigma) $, and write
$ G=\Gamma(\omega) $.
We suppose that $ G $ is not two-ended.
By Proposition 4.3, $ G $ is quasiconvex, and hence intrinsically
hyperbolic.
We shall, in turn, want to consider splittings of $ G $ over
finite groups, so to avoid any confusion later on, we shall
alter our notation, so that it is specific to this situation.
\par
Let $ \Xi $ be an indexing set which is in bijective correspondence
with the set, $ {\vec \Delta}(\omega) $, of directed edges of
$ \Sigma $ with heads at $ \omega $.
Thus, $ G $ permutes the elements of $ \Xi $.
There are finitely many $ G $--orbits (since
$ {\vec \Delta}(\omega)/\Gamma(\omega) $ is finite).
Given $ \xi \in \Xi $, we write $ H(\xi) $ for the stabiliser, in
$ G $, of $ \xi $.
Thus, if $ {\vec e} \in {\vec \Delta}(\omega) $ is the edge corresponding
to $ \xi $, then $ H(\xi)=\Gamma(e) $.
In particular, $ H(\xi) $ is two-ended.
Let $ J(\xi) = \Psi({\vec e}) $.
Thus, $ J(\xi) $ is a closed $ H(\xi) $--invariant subset of
$ \Lambda G $.
Moreover, $ \fr J(\xi) = J(\xi) \cap \Lambda G =\Lambda H(\xi) $ consists
of a pair of distinct points.
\par
In this notation, we have:
\gap
\lemma{6.1}
{\sl
$ \partial \Gamma = \Lambda G \cup \bigcup_{\xi \in \Xi} J(\xi) $.
}
\endprf
\gap
\lemma{6.2}
{\sl
$ \diam J(\xi) \rightarrow 0 $ for $ \xi \in \Xi $.
}
\endprf
\par
Here, Lemma 6.1 is a rewriting of Lemma 4.4, and Lemma 6.2 is
a restriction of Lemma 4.5.
\par
If $ K \subseteq \Lambda G $ is closed, we write
$ \Xi(K) = \{ \xi \in \Xi \mid \fr J(\xi) \subseteq K \} $, and write
$ \Upsilon(K) = K \cup \bigcup_{\xi \in \Xi(K)} J(\xi) $.
Thus, Lemma 4.6 says that:
\gap
\lemma{6.3}
{\sl
$ \Upsilon(K) $ is a closed subset of $ \partial \Gamma $.
}
\endprf
\par
These observations tell us all we need to know about the groups
$ H(\xi) $ and sets $ J(\xi) $ for the rest of this section.
Thus, for the moment, we can forget how they were constructed.
\par
Now, $ G $ is intrinsically hyperbolic, with $ \partial G $
identified with $ \Lambda G $.
We write $ \Lambda G = \Lambda_0 G \sqcup \Lambda_\infty G $,
corresponding to the partition
$ \partial G = \partial_0 G \sqcup \partial_\infty G $, as described
in Section 5.
Let $ T $ be a complete $ G $--tree, so that $ \partial T \equiv
\Lambda_\infty G $.
We write $ V_{\rm f\/in}(T) $ and $ V_{\rm inf}(T) $ respectively,
for the sets of vertices of $ T $ of finite and infinite degree.
Thus, $ \Lambda_0 G = \bigsqcup_{v \in V(T)} \Lambda G(v) $.
We note that if $ T $ is non-trivial (ie\ not a point), then
$ \Lambda_\infty G $ is dense in $ \Lambda G $.
\par
Given $ \xi \in \Xi $, the subgroup $ H(\xi) $ is two-ended.
It is either elliptic or hyperbolic with respect to the $ G $--tree
$ T $.
We write $ \Xi_{\rm ell} $ and $ \Xi_{\rm hyp} $, respectively,
for the sets of $ \xi \in \Xi $ such that $ H(\xi) $ is elliptic or
hyperbolic.
\par
If $ \xi \in \Xi_{\rm ell} $, then $ H(\xi) $ fixes a unique vertex
$ v(\xi) \in V_{\rm inf}(T) $, so that $ H(\xi) \subseteq G(v(\xi)) $
and $ \fr J(\xi) \subseteq \Lambda G(v(\xi)) $.
Given $ v \in V(T) $, we write $ \Xi_{\rm ell}(v) =
\{ \xi \in \Xi \mid H(\xi) \subseteq G(v) \} $.
Thus $ \Xi_{\rm ell}(v) \subseteq \Xi_{\rm ell} $, and
$ \Xi_{\rm ell}(v) = \emptyset $ for all $ v \in V_{\rm f\/in}(T) $.
In fact, $ \Xi_{\rm ell} = \bigsqcup_{v \in V(T)} \Xi_{\rm ell}(v) $.
\par
Given $ \xi \in \Xi_{\rm hyp} $, we write $ \beta(\xi) \subseteq T $ for
the unique biinfinite arc in $ T $ preserved setwise by $ H(\xi) $.
Note that, under the identification of $ \partial T $ and
$ \Lambda_0 G $, we have $ \partial \beta(\xi) = \Lambda H(\xi) $.
\par
Suppose that $ F \subseteq T $ is a finitely separated subgraph.
Recall from Section 5 that $ A(F) $ is defined as
$ A(F) = \partial F \cup \bigcup_{v \in V(F)} \Lambda G(v) $.
Thus, by Lemma 5.3, $ A(F) $ is closed in $ \Lambda G $ and
hence in $ \partial \Gamma $.
We abbreviate $ A(\Phi({\vec e})) $ to $ A({\vec e}) $.
(So that $ A({\vec e}) $ has the form $ \Psi({\vec e}) $ in the notation of
Section 4.)
\par
If $ F \subseteq T $ is finitely separated, we write
$ \Xi(F) = \Xi(A(F)) = \{ \xi \in \Xi \mid \fr J(\xi) \subseteq A(F) \} $.
Thus, $ \xi \in \Xi_{\rm ell} \cap \Xi(F) $ if and only if $ v(\xi) \in
V(F) $.
Also, $ \xi \in \Xi_{\rm hyp} \cap \Xi(F) $ if and only if
$ \partial \beta(\xi) \subseteq \partial F $.
\par
If $ {\vec e} \in {\vec E}(T) $, we shall abbreviate
$ \Xi({\vec e}) = \Xi(\Phi({\vec e})) $.
Thus, $ \xi \in \Xi({\vec e}) $ if and only if $ {\vec e} $ points
away from $ v(\xi) $ or $ \beta(\xi) $.
Suppose $ v_0 \in V(T) $.
Let $ \alpha \subseteq T $ be the arc joining $ v_0 $ to $ v(\xi) $
or to the nearest point of $ \beta(\xi) $.
Then, $ \{ {\vec e} \in {\vec \Omega}(v_0) \mid
\xi \in \Xi({\vec e}) \} $ consists of the directed edges in $ \alpha $
which point towards $ v_0 $.
In particular, this set is finite.
Indeed, if $ \Xi_0 \subseteq \Xi $ is finite, we see that
$ \{ {\vec e} \in {\vec \Omega}(v_0) \mid \Xi_0 \cap \Xi({\vec e})
\ne \emptyset \} $ is finite.
\par
If $ F \subseteq T $ is a finitely separated subgraph, we write
$$ B(F) = A(F) \cup \bigcup_{\xi \in \Xi(F)} J(\xi) .$$
In other words, $ B(F) = \Upsilon(A(F)) $, as defined earlier in
this section.
Thus, by Lemma 6.3, we have:
\gap
\lemma{6.4}
{\sl
The set $ B(F) \subseteq \partial \Gamma $ is closed, for any
finitely separated subgraph, $ F $, of $ T $.
}
\endprf
\par
If $ {\vec e} \in {\vec E}(T) $, we abbreviate $ B({\vec e})
= B(\Phi({\vec e})) $.
\gap
\lemma{6.5}
{\sl
If $ v_0 \in V(T) $, then $ \diam B({\vec e}) \rightarrow 0 $
for $ {\vec e} \in {\vec \Omega}(v_0) $.
}
\gap
\prf
Suppose $ \delta > 0 $.
By Lemma 6.2, there is a finite subset $ \Xi_0 \subseteq \Xi $
such that if $ \xi \in \Xi \bksl \Xi_0 $ then
$ \diam J(\xi) \le \delta/3 $.
Let $ {\vec \Omega}_0 = \{ {\vec e} \in {\vec \Omega}(v_0) \mid
\Xi_0 \cap \Xi({\vec e}) \ne \emptyset \} $.
As observed above, $ {\vec \Omega}_0 $ is finite.
Let $ {\vec \Omega}_1 = \{ {\vec e} \in {\vec \Omega}(v_0) \mid
\diam A({\vec e}) \ge \delta/3 \} $.
By Lemma 4.5, $ {\vec \Omega}_1 $ is also finite.
\par
Suppose $ {\vec e} \in {\vec \Omega}(v_0) \bksl
({\vec \Omega}_0 \cup {\vec \Omega}_1) $.
Suppose $ x \in B({\vec e}) $.
If $ x \notin A({\vec e}) $, then $ x \in J(\xi) $ for some
$ \xi \in \Xi({\vec e}) $.
Since $ {\vec e} \notin {\vec \Omega}_0 $, $ \Xi_0 \cap \Xi({\vec e}) =
\emptyset $, so $ \xi \notin \Xi_0 $.
Therefore, $ \diam J(\xi) \le \delta/3 $.
Now, $ \fr J(\xi) \subseteq A({\vec e}) $, and so
$ \rho(x,A({\vec e})) \le \delta/3 $.
This shows that $ B({\vec e}) $ lies in a $ (\delta/3) $--neighbourhood
of $ A({\vec e}) $.
Now, since $ {\vec e} \notin {\vec \Omega}_1 $, $ \diam A({\vec e})
< \delta/3 $ and so $ \diam B(\epsilon) < \delta $.
\endprf
\par
Recall, from Section 3, that if $ S \subseteq T $ is a subtree,
then there is a natural projection $ \pi_S \co T \cup \partial T
\longrightarrow S \cup \partial S $.
If $ v \in V(S) $, we write $ F(S,v) $ for the subtree
$ T \cap \pi_S^{-1} v $.
If $ R \subseteq S $ is a subtree, then we see that
$ F(S,v) \subseteq F(R, \pi_R v) $.
Recall that $ {\vec \Delta}(S) = \{ {\vec e} \in {\vec E}(T) \mid
\rhead({\vec e}) \in S, \rtail({\vec e}) \notin S \} $.
If $ v \in V(S) $, set $ {\vec \Delta}(S,v) = {\vec \Delta}(S) \cap
{\vec \Delta}(v) $.
We write $ {\vec \Omega}(S) $ for the set of all directed edges pointing
towards $ S $, ie\ $ {\vec \Omega}(S) = \bigcap_{v \in V(S)}
{\vec \Omega}(v) $.
Clearly, $ {\vec \Delta}(S) \subseteq {\vec \Omega}(S) $.
Also if $ R \subseteq S $ is a subtree, then
$ {\vec \Omega}(S) \subseteq {\vec \Omega}(R) $.
If $ v \in V(T) \bksl V(R) $, let $ {\vec e}(R,v) $ be the directed
edge with head at $ \pi_R v $ which lies in the arc joining $ v $
to $ \pi_R v $.
In other words, $ {\vec e}(R,v) $ is the unique edge in
$ {\vec \Delta}(R) $ such that $ v \in \Phi({\vec e}(R,v)) $.
Note that, if $ v \in V(S) \bksl V(R) $, then $ F(S,v) \subseteq
\Phi({\vec e}(R,v)) $.
\par
Let $ {\script T} $ be the set of all finite subtrees of $ T $.
Given $ \delta > 0 $, let
$$
\eqalign{
{\script T}_1(\delta) &= \{ S \in {\script T} \mid (\forall {\vec e} \in
{\vec \Delta}(S))(\diam B({\vec e}) < \delta) \} \cr
{\script T}_2(\delta) &= \{ S \in {\script T} \mid (\forall v \in
V(S) \cap V_{\rm f\/in}(T))(\diam B(F(S,v)) < \delta) \} \cr
{\script T}_3(\delta) &= \{ S \in {\script T} \mid (\forall v \in
V(S) \cap V_{\rm inf}(T))(\forall {\vec e} \in {\vec \Delta}(S,v))
(\rho(\Lambda G(v), B({\vec e})) < \delta) \} .}
$$
Let $ {\script T}(\delta) = {\script T}_1(\delta) \cap
{\script T}_2(\delta) \cap {\script T}_3(\delta) $.
\par
It is really the collection $ {\script T}(\delta) $ in which we
shall ultimately be interested.
It can be described a little more directly as follows.
A finite tree, $ S $, lies in $ {\script T}(\delta) $ if and only if for
each $ v \in V(S) $, we have either $ v \in V_{\rm f\/in}(T) $
and $ \diam B(F(S,v)) < \delta $ or else $ v \in V_{\rm inf}(T) $
and for all $ {\vec e} \in {\vec \Delta}(S,v) $ we have
$ \diam B({\vec e}) < \delta $ and
$ \rho(\Lambda G(v), B({\vec e})) < \delta $.
It is this formulation we shall use in applications.
\par
Note that if $ R \in {\script T}_1(\delta) $, then, in fact,
$ \diam B({\vec e}) < \delta $ for all $ {\vec e} \in {\vec \Omega}(R) $.
We see that if $ R \in {\script T}_1(\delta) $, $ S \in {\script T} $
and $ R \subseteq S $, then $ S \in {\script T}_1(\delta) $.
More to the point, we have:
\gap
\lemma{6.6}
{\sl
If $ R \in {\script T}(\delta) $, $ S \in {\script T} $ and
$ R \subseteq S $, then $ S \in {\script T}(\delta) $.
}
\gap
\prf
As observed above, $ S \in {\script T}_1(\delta) $.
\par
Suppose that $ v \in V(S) \cap V_{\rm f\/in}(T) $.
If $ v \in V(R) $, then $ F(S,v) \subseteq F(R,v) $, and so
$ B(F(S,v)) \subseteq B(F(R,v)) $.
Therefore, $ \diam B(F(S,v)) \le \diam B(F(R,v)) < \delta $, since
$ R \in {\script T}_2(\delta) $.
On the other hand, if $ v \notin V(R) $, then $ F(S,v) \subseteq
\Phi({\vec e}(R,v)) $, so $ \diam B(F(S,v)) \le
\diam B({\vec e}(R,v)) < \delta $, since $ R \in {\script T}_1(\delta) $.
This shows that $ S \in {\script T}_2(\delta) $.
\par
Finally, suppose $ v \in V(S) \cap V_{\rm inf}(T) $ and
$ {\vec e} \in {\vec \Delta}(S,v) $.
If $ v \in V(R) $, then $ {\vec e} \in {\vec \Delta}(R,v) $, so
$ \rho(\Lambda G(v),B({\vec e})) $, since $ R \in
{\script T}_3(\delta) $.
On the other hand, if $ v \notin V(R) $, then $ \{ v \} \cup
\Phi({\vec e}) \subseteq F(R,{\vec e}(R,v)) $, and so
$ \Lambda G(v) \cup B({\vec e}) \subseteq B(F(R,{\vec e}(R,v))) $.
But $ \diam B(F(R,{\vec e}(R,v))) < \delta $, since $ R \in
{\script T}_1(\delta) $.
In particular, $ \rho(\Lambda G(v), B({\vec e})) < \delta $.
This shows that $ S \in {\script T}_3(\delta) $.
\endprf
\gap
\lemma{6.7}
{\sl
$ {\script T}(\delta) \ne \emptyset $.
}
\gap
\prf
Using Lemma 6.5, we can certainly find some $ R \in
{\script T}_1(\delta) $.
We form another finite tree, $ S \supseteq R $, by adjoining a finite
number of adjacent edges as follows.
If $ v \in V(R) \cap V_{\rm f\/in}(T) $, we add all edges which are
incident on $ v $.
If $ v \in V(R) \cap V_{\rm inf}(T) $, we add all those incident edges,
$ e $, which correspond to $ {\vec e} \in {\vec \Delta}(R,v) $ for which
$ \rho(\Lambda G(v),B({\vec e})) \ge \delta $.
By Lemma 4.1, and the fact that $ {\vec \Delta}(v)/G(v) $ is finite,
there are only finitely many such $ {\vec e} $.
We thus see that $ S $ is finite.
The fact that $ S \in {\script T}(\delta) $ follows by essentially the
same arguments as were used in the proof of Lemma 6.6.
\endprf
\gap\gap
\section{Connectedness properties of boundaries of hyperbolic groups}
\gap\par
In this section, we continue the analysis of Section 6, bringing
connectedness assumptions into play.
\par
Suppose, as before, that $ \Gamma $ is a hyperbolic group, and
that $ \Sigma $ is a cofinite $ \Gamma $--tree with two-ended
edge stabilisers.
We now add the assumption that $ \Gamma $ is one ended, so that
$ \partial \Gamma $ is a continuum.
In this case, we note:
\gap
\lemma{7.1}
{\sl
For each $ {\vec e} \in {\vec E}(\Sigma) $, the set $ \Psi({\vec e}) $
is connected.
}
\gap
\prf
Since $ \Gamma(e) $ is two-ended, we have $ \fr \Psi({\vec e}) =
\Lambda \Gamma(e) = \{ a,b \} $, where $ a,b \in \Psi({\vec e}) $
are distinct.
Moreover, $ \Psi({\vec e}) $ is closed and $ \Gamma(e) $--invariant.
Also $ \Psi({\vec e}) \ne \{ a,b \} $, since it must, for example,
contain all points of $ \partial \Phi({\vec e}) $.
\par
Let $ K $ be a connected component of $ \Psi({\vec e}) $.
We claim that $ K \cap \{ a,b \} \ne \emptyset $.
To see this, suppose $ a,b \notin K $.
There are subsets $ K_1, K_2 \subseteq \Psi({\vec e}) $,
containing $ K $, with $ a \notin K_1 $, $ b \notin K_2 $, and which
are clopen in $ \Psi({\vec e}) $.
Let $ L = K_1 \cap K_2 $.
Thus, $ K \subseteq L \subseteq \Psi({\vec e}) \bksl \fr
\Psi({\vec e}) $.
Since $ \Psi({\vec e}) $ is closed in $ \partial \Gamma $,
so is $ L $, and since $ \Psi({\vec e}) \bksl \partial \Psi({\vec e}) $
is open in $ \partial \Gamma $, so also is $ L $.
In other words, $ L $ is clopen in $ M $, contradicting the
hypothesis that $ \partial \Gamma $ is connected.
\par
Suppose, then, that $ a \in K $.
Let $ H \le \Gamma(e) $ be the subgroup (of index at most 2) fixing
$ a $ (and hence $ b $).
We see that $ K $ is $ H $--invariant.
Now $ \Lambda H = \{ a,b \} $ so either $ b \in K $, or $ K = \{ a \} $.
In the former case, we see that $ K = \Psi({\vec e}) $, showing
that $ \Psi({\vec e}) $ is connected.
In the latter case, we see, by a similar argument, that
the component of $ K $ containing $ b $ equals $ \{ b \} $, giving
the contradiction that $ \Psi({\vec e}) = \{ a,b \} $.
\endprf
\par
Now, as in Section 6, we focus on one vertex $ \omega \in V(\Sigma) $,
and write $ G = \Gamma(\omega) $.
Let $ T $ be a complete $ G $--tree.
Now, $ \Lambda G = \Lambda_0 G \sqcup \Lambda_\infty G $, where
$ \Lambda_0 G = \bigsqcup_{v \in V(T)} \Lambda G(v) $ and
$ \Lambda_\infty G $ is identified with $ \partial T $.
It is possible that $ T $ may be trivial, but most of the following
discussion will be vacuous in that case.
If not, then $ \Lambda_\infty G $ is dense in $ \Lambda G $.
\par
We now reintroduce the notation used in Section 6, namely
$ \Xi $, $ J(\xi) $, $ H(\xi) $, $ B({\vec e}) $, etc.
Note that if $ \xi \in \Xi $ corresponds to the directed edge
$ {\vec \epsilon} $ of $ \Sigma $, then $ J(\xi) $ equals
$ \Psi({\vec \epsilon}) $ and the closure of
$ \partial \Gamma \bksl J(\xi) $ in $ \partial \Gamma $ equals
$ \Psi(-{\vec \epsilon}) $ (in the notation of Section 4).
Thus, rephrasing Lemma 7.1, we get:
\gap
\lemma{7.2}
{\sl
For each $ \xi \in \Xi $, the set $ J(\xi) $ is connected.
Moreover, the closure of $ \partial \Gamma \bksl J(\xi) $
in $ \partial \Gamma $ is also connected.
}
\endprf
\par
Let $ {\script B} = \{ \beta(\xi) \mid \xi \in \Xi_{\rm hyp} \} $.
Now, $ \Xi_{\rm hyp}/G $ is finite, so Lemma 2.1 tells us that:
\gap
\lemma{7.3}
{\sl
The arc system $ {\script B} $ is edge-finite.
}
\endprf
\par
Now, since $ \Gamma $ is one-ended, the set of two-ended subgroups
$ {\script H} = \{ H(\xi) \mid \xi \in \Xi_{\rm hyp} \} $ is
indecomposable.
Since $ {\script B} $ is the set of axes of elements of $ {\script H} $,
we see by Proposition 3.8 that:
\gap
\lemma{7.4}
{\sl
$ {\script B} $ is indecomposable.
}
\endprf
\par
Alternatively, one can give a direct proof of Lemma 7.4 along the lines
of Lemma 5.4.
Thus, if $ {\script B} $ is decomposable, we can find two finitely
separated subgraphs, $ F_1 $ and $ F_2 $, of $ T $, so that
$ \partial G = A(F_1) \sqcup A(F_2) $, and such that for all
$ \xi \in \Xi_{\rm hyp} $, either $ \partial \beta(\xi) \in \partial F_1 $,
or $ \partial \beta(\xi) \in \partial F_2 $.
It follows that $ \partial \Gamma = B(F_1) \sqcup B(F_2) $ are closed
in $ \partial \Gamma $, contradicting the assumption that
$ \partial \Gamma $ is connected.
\gap
\par
To go further, we shall want some more general observations and notation
regarding simplicial trees.
For the moment, $ T $ can be any simplicial tree, and $ {\script B} $ any
arc system on $ T $.
\par
In Section 3, we associated to any finite subtree, $ S \subseteq T $,
an equivalence relation, $ \mathord{\sim}_S = \mathord{\sim}_{S,{\script B}} $,
on $ \partial T $.
This, in turn, gives us a subpartition, $ {\script W}(S) $, of the set
$ V(S) $ of vertices of $ S $.
The elements of $ {\script W}(S) $ are the vertex sets of the
connected components of the Whitehead graph, $ {\script G}(S) $.
\par
More generally, we shall say that a subtree, $ S $, of $ T $ is
\em{bounded} if it has finite diameter in the combinatorial metric
on $ T $.
In particular, every arc of $ {\script B} $ meets $ S $, if at all,
in a compact interval (or point).
We define the equivalence relation, $ \mathord{\sim}_S =
\mathord{\sim}_{S,{\script B}} $ on $ \partial T $ in exactly the
same way as for finite trees.
We also get a graph $ {\script G}(S) $, and a subpartition,
$ {\script W}(S) $ of $ V(S) $ as before.
Note that if $ {\script B} $ is edge-finite, then
$ {\script G}(S) $ is locally finite.
\par
We have already observed that if $ R \subseteq S $ is a subtree of $ S $,
then the relation $ \mathord{\sim}_R $ is coarser than
the relation $ \mathord{\sim}_S $ (ie\ $ x \sim_S y $ implies
$ x \sim_R y $).
Moreover, the subpartition, $ {\script W}(R) $ of $ V(R) $ can
be described explicitly in terms of the subpartition $ {\script W}(S) $
and the map $ \pi_R|V(S) \co V(S) \longrightarrow V(R) $.
To do this, define $ \mathord{\cong} $ to be the equivalence relation
on $ {\script W}(S) $ generated by relations of the form
$ W \cong W' $ whenever $ \pi_R W \cap \pi_R W' \ne \emptyset $.
An element of $ {\script W}(R) $ is then a union of sets of the
form $ \pi_R W $ as $ W $ ranges over some $ \mathord{\cong} $--class
in $ {\script W}(S) $.
For future reference, we note:
\gap
\lemma{7.5}
{\sl
Suppose $ R \subseteq S $ are bounded subtrees of $ T $.
If $ W \in {\script W}(S) $, $ W \subseteq V(R) $, and
$ W \cap \pi_R(V(S) \bksl V(R)) = \emptyset $, then $ W \in {\script W}(R) $.
}
\gap
\prf
If $ W' \in {\script W}(S) $ and $ W \cap \pi_R W' \ne \emptyset $,
then $ W \cap W' \ne \emptyset $.
(To see this, choose $ v \in W' $ with $ \pi_R v \in W \subseteq V(R) $.
Since $ W \cap \pi_R(V(S) \bksl V(R)) = \emptyset $, it follows
that $ v \in V(R) $, so $ \pi_R v = v $.
Thus $ v \in W \cap W' $.)
Since $ W,W' \in {\script W}(S) $ we thus have $ W=W' $,
so $ W' = \pi_R W' $.
This shows that any set of the form $ \pi_R W' $ for $ W' \in
{\script W}(S) $ which meets $ W $ must, in fact, be equal to $ W $.
 From the description of $ {\script W}(R) $ given above, we see
that $ W \in {\script W}(R) $.
\endprf
\par
Given a directed edge $ {\vec e} \in {\vec E}(T) $, let
$ {\script S}({\vec e}) $ be the set of finite subtrees, $ S $, of
$ T $ with the property that $ {\vec \Delta}(\rhead({\vec e}))
\cap {\vec E}(S) = \{ {\vec e} \} $
(ie\ $ e \subseteq S $, and $ \rhead({\vec e}) $ is a terminal
vertex of $ S $).
Given $ S \in {\script S}({\vec e}) $, we define the equivalence
relation $ \mathord{\simeq}_S $ on $ \partial \Phi({\vec e}) $
to be the transitive closure of relations of the form
$ x \simeq_S y $ whenever $ \pi_S x = \pi_S y $ or $ \partial \beta
= \{ x,y \} $ for some $ \beta \in {\script B} $, with $ \beta
\subseteq \Phi({\vec e}) $.
Clearly, if $ x \simeq_S y $ then $ x \sim_S y $.
Also, if $ R,S \in {\script S}({\vec e}) $ with $ R \subseteq S $,
then $ x \simeq_S y $ implies $ x \simeq_R y $.
We can also define a subpartition, $ {\script W}(S,{\vec e}) $, of
$ V(S) \bksl \{ \rhead({\vec e}) \} $, in a similar manner to
$ {\script W}(S) $, as described in Section 3.
\par
Suppose now that $ {\script B} $ is edge-finite and indecomposable,
and suppose $ S \in {\script S}({\vec e}) $.
Suppose $ Q \subseteq \partial \Phi({\vec e}) $ is a
$ \mathord{\simeq}_S $--class.
Since there is only one $ \mathord{\sim}_S $--class, there must
be some $ \beta \in {\script B} $ with one endpoint in $ Q $ and
one endpoint in $ \partial \Phi(-{\vec e}) $.
Thus, $ e \subseteq \beta $.
It follows that the number of $ \mathord{\simeq}_S $--classes is
bounded by the number of arcs in $ {\script B} $ containing the edge
$ e $.
By the edge-finiteness assumption, this number is finite.
It follows that, as the trees $ S \in {\script S}({\vec e}) $ get
bigger, the relations $ \mathord{\simeq}_S $ must stabilise.
More precisely, there is a (unique) equivalence relation,
$ \mathord{\simeq} $, on $ \partial \Phi({\vec e}) $ such that
the set $ {\script S}_0({\vec e}) = \{ S \in {\script S}({\vec e}) \mid
\mathord{\simeq}_S = \mathord{\simeq} \} $ contains all but finitely
many elements of $ {\script S}({\vec e}) $.
Note that if $ R \in {\script S}_0({\vec e}) $, $ S \in
{\script S}({\vec e}) $, and $ R \subseteq S $, then
$ S \in {\script S}_0({\vec e}) $.
Note also that there are finitely many $ \mathord{\simeq} $--classes.
\gap
\par
We now return to the set-up described earlier, with $ T $ a complete
$ G $--tree, and with $ {\script B} = \{ \beta(\xi) \mid
\xi \in \Xi_{\rm hyp} \} $.
We have seen that $ {\script B} $ is edge-finite and indecomposable.
We note:
\gap
\lemma{7.6}
{\sl
Suppose $ {\vec e} \in {\vec E}(T) $ and $ x,y \in
\partial \Phi({\vec e}) $.
If $ x \simeq y $, then $ x $ and $ y $ lie in the same connected
component of $ B({\vec e}) $.
}
\gap
\prf
Suppose, for contradiction that $ x $ and $ y $ lie in different
components of $ B({\vec e}) $.
We can partition $ B({\vec e}) $ into two closed
subsets, $ B({\vec e}) = K \sqcup L $, with $ x \in K $ and $ y \in L $.
\par
Let $ \delta = {1 \over 2} \rho(K,L) >0 $.
By Lemma 6.7, we can find some $ R \in {\script T}(\delta) $.
By Lemma 6.6, we can suppose that $ S = R \cap (e \cup \Phi({\vec e}))
\in {\script S}_0({\vec e}) $.
(For example, take $ R $ to be the smallest tree containing a given
element of $ {\script T}(S) $ and a given element of
$ {\script S}_0({\vec e}) $.)
Thus, $ \mathord{\simeq}_S = \mathord{\simeq} $, so in particular,
$ x \simeq_S y $.
Note that, if $ v \in V(S) \bksl \{ \rhead({\vec e}) \} $, then
$ F(R,v) = F(S,v) $ (in the notation of Section 2).
\par
Now, from the definition of the relation $ \mathord{\simeq}_S $, we
have a finite sequence, $ x=x_0, x_1 ,\ldots, x_n =y $ of points
of $ \partial \Phi({\vec e}) $, such that for each $ i $, either
$ \pi_S x_i = \pi_S x_{i+1} $, or there is some $ \xi \in \Xi_{\rm hyp} $,
with $ \partial \beta(\xi) = \{ x_i, x_{i+1} \} $.
Now, $ \partial \Phi({\vec e}) \subseteq B({\vec e}) = K \sqcup L $,
so for each $ i $, either $ x_i \in K $ or $ x_i \in L $.
We claim, by induction on $ i $, that $ x_i \in K $ for all $ i $.
\par
Suppose, then, that $ x_i \in K $.
Suppose first, that $ \{ x_i, x_{i+1} \} = \partial \beta(\xi) $ for
some $ \xi \in \Xi_{\rm hyp} $.
We have that $ x_i, x_{i+1} \in J(\xi) \subseteq B({\vec e}) $.
Moreover, by Lemma 6.1, $ J(\xi) $ is connected.
It follows that $ x_{i+1} \in K $.
\par
We can thus suppose that $ \pi_S x_i = \pi_S x_{i+1} = v \in
V(S) \bksl \{ \rhead({\vec e}) \} $.
Thus, $ x_i, x_{i+1} \in \partial F(S,v) = \partial F(R,v) \subseteq
B(F(R,v)) $.
Now, if $ v \in V_{\rm f\/in}(T) $, then, since $ R \in
{\script T}(\delta) $, we have $ \diam B(F(R,v)) < \delta $.
Therefore, $ \rho(x_i, x_{i+1}) < \delta $ and so $ x_{i+1} \in K $.
Thus, we can assume that $ v \in V_{\rm inf}(T) $.
Since $ x_i \in \partial F(R,v) $, we have $ x_i \in
\partial \Phi({\vec \epsilon}) $ for some $ {\vec \epsilon} \in
{\vec \Delta}(R,v) $.
Again, since $ R \in {\script T}(\delta) $, we have
$ \diam B({\vec \epsilon}) < \delta $ and
$ \rho(B({\vec \epsilon}), \Lambda G(v)) < \delta $.
Thus, $ \rho(x_i,\Lambda G(v)) < 2 \delta $.
Similarly, $ \rho(x_{i+1}, \Lambda G(v)) < 2 \delta $.
Now, $ \Lambda G(v) $ is connected, and so it again follows that
$ x_{i+1} \in K $.
\par
Thus, by induction on $ i $, we arrive at the contradiction that
$ y=x_n \in K $.
This shows that $ x $ and $ y $ lie in the same component of
$ B({\vec e}) $ as required.
\endprf
\par
Now, fix some $ v \in V_{\rm inf}(T) $, so that $ G(v) $ is
one-ended, and $ \Lambda G(v) $ is a subcontinuum of
$ \partial \Gamma $.
\par
We say that a $ G(v) $--invariant subtree, $ S $, of $ T $ is
{\it stable about\/} $ v $ if $ S \cap \Phi({\vec e}) \in
{\script S}_0({\vec e}) $ for all $ {\vec e} \in {\vec \Delta}(v) $.
Note that, since $ {\vec \Delta}(v)/G(v) $ is finite,
$ S/G(v) $ is finite.
In particular, we see that $ S $ is bounded (ie\ has finite
diameter).
Note that, since $ S $ contains every edge of $ T $ incident on
$ v $, we have $ \pi_S \partial T \subseteq V(S) \bksl \{ v \} $.
Let $ \mathord{\sim}_S = \mathord{\sim}_{S,{\script B}} $ be the equivalence
relation on $ \partial T $ as defined in Section 3 (in the case of finite
trees).
We remark that $ \mathord{\sim}_S $ is independent of the
choice of stable tree, $ S $, since it is easily seen to be
definable purely in terms of the arc system $ {\script B} $, and
the relations, $ \mathord{\simeq} $ for $ {\vec e} \in
{\vec \Delta}(v) $.
We shall thus write $ \mathord{\sim}_S $ simply as $ \mathord{\sim} $.
Clearly, $ \mathord{\sim} $ is $ G(v) $--invariant.
(It need not be trivial, since we are only assuming that
$ S $ is bounded.)
\par
We can certainly construct a stable tree about $ v $ by
taking $ S = \bigcup_{{\vec e} \in {\vec \Delta}(v)} S({\vec e}) $.
In this case, $ S \cap \Phi({\vec e}) = S({\vec e}) \in
{\script S}_0({\vec e}) $.
\par
Note that we get a subpartition, $ {\script W}(S) $, of
$ V(S) $, as described in Section 3.
Note that $ \bigcup {\script W}(S) \subseteq \pi_S \partial T $.
In particular, $ v \notin \bigcup {\script W}(S) $.
\gap
\lemma{7.7}
{\sl
The setwise stabiliser, in $ G(v) $, of every $ \mathord{\sim} $--class
is infinite.
}
\gap
\prf
As described in Section 3, each $ \mathord{\sim} $--class corresponds
to an element of $ {\script W}(S) $.
Moreover, $ (\bigcup {\script W}(S))/G(v) \subseteq V(S)/G(v) $ is
finite.
Thus, the lemma is equivalent to asserting that each element of
$ {\script W}(S) $ is infinite.
\par
Suppose, to the contrary, that $ W \in {\script W}(S) $ is finite.
Let $ {\vec \Delta}_0 = \{ {\vec e} \in {\vec \Delta}(v) \mid
W \cap S({\vec e}) \ne \emptyset \} $, and let
$ R = \bigcup_{{\vec e} \in {\vec \Delta}_0} S({\vec e}) $.
Thus, $ R $ is a finite subtree of $ S $, and $ W \subseteq V(R) $.
Moreover, $ \pi_R(V(S) \bksl V(R)) = \{ v \} $, so, in particular,
$ W \cap \pi_R(V(S) \bksl V(R)) = \empty $.
Thus, by Lemma 7.5, $ W \in {\script W}(R) $.
But $ v \in \bigcup {\script W}(R) $ (since any element of
$ \partial \Phi({\vec e}) $ for $ {\vec e} \in {\vec \Delta}(v) \bksl
{\vec \Delta}_0 $ projects to $ v $ under $ \pi_R $).
Thus, $ {\script W}(R) \ne \{ W \} $.
This shows that there is more than one $ \mathord{\sim}_R $--class,
contradicting the fact that $ {\script B} $ is indecomposable.
\endprf
\par
Finally, we note:
\gap
\lemma{7.8}
{\sl
If $ x,y \in \partial T $ with $ x \sim y $, then $ x $ and $ y $
lie in the same quasicomponent of $ \partial \Gamma \bksl \Lambda G(v) $.
}
\gap
\prf
In fact, we shall show that $ x $ and $ y $ both lie in a compact connected
subset, $ K $, of $ \partial \Gamma \bksl \Lambda G(v) $.
\par
By the definition of the relation $ \mathord{\sim} = \mathord{\sim}_S $,
we can assume that either $ \pi_S x = \pi_S y $ or
there is some $ \xi \in \Xi_{\rm hyp} $ with
$ \partial \beta(\xi) = \{ x,y \} $.
\par
In the former case, let $ w = \pi_S x = \pi_S y $.
Thus, $ w \in V(S({\vec e})) $ for some $ {\vec e} \in {\vec \Delta}(v) $.
Since $ S({\vec e}) \in {\script S}_0({\vec e}) $, we have
$ x \simeq y $, and so, by Lemma 7.6, $ x $ and $ y $ lie in the
same component of $ B({\vec e}) $.
Call this component $ K $.
Thus, $ K $ is closed in $ B({\vec e}) $ and hence in $ \partial \Gamma $.
Note that, from the definition of $ B({\vec e}) $, we have
$ B({\vec e}) \cap \Lambda G(v) = \emptyset $ and so
$ K \cap \Lambda G(v) = \emptyset $.
\par
In the latter case, set $ K = J(\xi) $.
Thus, by Lemma 6.1, $ K $ is connected.
Also $ K \cap \Lambda G = \{ x,y \} \subseteq \partial T $, and
so, again, $ K \cap \Lambda G(v) = \emptyset $.
\endprf
\gap\gap
\section{Global cut points}
\gap\par
In this section, we set out the ``inductive step'' of the proof
that a strongly accessible hyperbolic group has no global cut
points in its boundary.
In the light of the result announced in [\DeP], we see that this, in
fact, applies to all one-ended hyperbolic groups.
A more direct proof of the general case was given in [\Swa]
using the results of [\Boa,\Boc,\L].
(See also [\Bod].)
\par
Specifically, we shall show:
\gap
\theorem{8.1}
{\sl
Suppose that $ \Gamma $ is a one-ended hyperbolic group.
Suppose that we represent $ \Gamma $ as a finite graph of groups over
two-ended subgroups.
Suppose that each maximal one-ended subgroup of each vertex group
has no global cut point in its boundary (as an intrinsic hyperbolic group).
Then, $ \partial \Gamma $ has no global cut point.
}
\gap
\ppar
Before we start on the proof, we give a few general definitions and
observations relating to global cut points.
\par
Suppose that $ M $ is any continuum, ie\ a compact connected
hausdorff space.
(For the moment, the compactness assumption is irrelevant.)
If $ p \in M $, and $ O,U \subseteq M $, we write $ OpU $ to
mean that $ O $ and $ U $ are non-empty open subsets and that
$ M $ is (set theoretically) a disjoint union $ M = O \sqcup \{ p \}
\sqcup U $.
Note that $ \fr O = \fr U = \{ p \} $.
Also, it's not hard to see that $ O \cup \{ p \} $ and
$ U \cup \{ p \} $ are connected.
(More discussion of this is given in [\Boa].)
We say that a point $ p \in M $ is a {\it global cut point\/} if
there exist $ O,U \subseteq M $ with $ OpU $.
\gap
\defne
If $ Q \subseteq M $ is any subset, and $ p \in M $, we say that
$ Q $ is {\it indivisible in $ M $ at $ p $ \/} if whenever we have
$ O,U \subseteq M $ with $ OpU $, then either $ Q \cap O = \emptyset $
or $ Q \cap U = \emptyset $.
\brk
If $ R \subseteq M $ is another subset, we say that $ Q $ is
{\it indivisible in $ M $ over $ R $ \/}, if it is indivisible in
$ M $ at every point of $ R $.
\brk
We say that $ Q $ is {\it (globally) indivisible\/} in $ M $ if
it is indivisible at every point of $ M $.
\gap
Thus, $ M $ is indivisible in itself if and only if it does not contain
a global cut point.
\par
Obviously, if $ P \subseteq Q \subseteq M $ and $ Q $ is indivisible in
$ M $, then so is $ P $.
Also any subcontinuum of $ M $ with no global cut point is indivisible in
$ M $.
We shall need the following simple observations:
\gap
\lemma{8.2}
{\sl
If $ P,Q \subseteq M $ are indivisible in $ M $, and
$ \card(P \cap Q) \ge 2 $, then $ P \cup Q $ is indivisible in $ M $.
}
\gap
\prf
Suppose $ OpU $.
Choose any $ x \in P \cap Q \bksl \{ p \} $.
We can assume that $ x \in O $, so that $ P \cap U = Q \cap U = \emptyset $.
Thus $ (P \cup Q) \cap U = \emptyset $.
\endprf
\gap
\lemma{8.3}
{\sl
Suppose that $ {\script Q} $ is a chain of indivisible subsets of
$ M $ (ie\ if $ P,Q \in {\script Q} $, then $ P \subseteq Q $ or
$ Q \subseteq P $).
Then $ \bigcup {\script Q} $ is indivisible.
}
\gap
\prf
Suppose $ OpU $, and $ x \in O \cap (\bigcup {\script Q}) $
and $ y \in U \cap (\bigcup {\script Q}) $.
Then $ x,y \in Q $ for some $ Q \in {\script Q} $, contradicting the
indivisibility of $ Q $.
\endprf
\gap
\lemma{8.4}
{\sl
If $ Q $ is indivisible in $ M $, then so is its closure, $ {\bar Q} $.
}
\gap
\prf
If $ OpU $, then we can assume that $ O \cap Q = \emptyset $, so
$ O \cap {\bar Q} = \emptyset $.
\endprf
\par
Now, let $ \Gamma $ be a one-ended hyperbolic group, and let
$ \Sigma $ be a cofinite $ \Gamma $--tree with two-ended edge
stabilisers.
We begin with the following observation:
\gap
\lemma{8.5}
{\sl
If $ \Lambda \Gamma(v) $ is indivisible in $ \partial \Gamma  $ for
all $ v \in V(\Sigma) $, then $ \partial \Gamma $ is indivisible.
}
\gap
\prf
Note that if $ v,w \in V(\Sigma) $ are adjacent, then
$ \Gamma(v) \cap \Gamma(w) $ is two-ended, so
$ \Lambda \Gamma(v) \cap \Lambda \Gamma(w) =
\Lambda(\Gamma(v) \cap \Gamma(w)) $ consists of a pair of points.
Thus, by Lemma 8.2, $ \Lambda \Gamma(v) \cap \Lambda \Gamma(w) $
is indivisible in $ \partial \Gamma $.
By an induction argument, we see that $ \bigcup_{v \in V(S)}
\Lambda \Gamma(v) $ is indivisible for any finite subtree,
$ S \subseteq \Sigma $.
Taking an exhaustion of $ \Sigma $ by an increasing sequence of finite
subtrees, and applying Lemma 8.3, we see that
$ \bigcup_{v \in V(\Sigma)} \Lambda \Gamma(v) $ is indivisible.
But this set is dense in $ \partial \Gamma $ (since it is non-empty
and $ \Gamma $--invariant).
The result follows by Lemma 8.4.\Endprf
\par
In fact, it's enough to verify the hypotheses of Lemma 8.5 for
those $ v \in V(\Sigma) $ for which $ \Gamma(v) $ is not two-ended.
To see this, first note that if $ \alpha $ is a finite arc connecting
two points $ v_0, v_1 \in V(\Sigma) $ such that $ \Gamma(v) $ is
two ended for all $ v \in V(\alpha) \bksl \{ v_0, v_1 \} $, then
the groups $ \Gamma(e) $ and $ \Gamma(v) $ are all commensurable for
all $ e \in E(\alpha) $ and $ v \in V(\alpha) \bksl \{ v_0, v_1 \} $.
Now, since $ \Gamma $ is hyperbolic and not two-ended, there must be
some $ v_0 \in V(\Sigma) $ such that $ \Gamma(v_0) $ is not
two-ended.
Suppose that $ v \in V(\Sigma) $ is some other vertex.
Connect $ v $ to $ v_0 $ by an arc in $ \Sigma $, and let $ w $ be
the first vertex of this arc for which $ \Gamma(w) $ is not two-ended.
Thus, $ \Gamma(v) \cap \Gamma(w) $ has finite index $ \Gamma(v) $,
and so $ \Lambda \Gamma(v) \subseteq \Lambda \Gamma(w) $.
Clearly, if $ \Lambda \Gamma(w) $ is indivisible in $ \partial \Gamma $,
then so is $ \Lambda \Gamma(v) $.
\par
As in Section 7, we now fix $ \omega \in V_{\rm inf}(\Sigma) $ and set
$ G = \Gamma(\omega) $.
We are interested in the indivisibility properties of $ \Lambda G $
as a subset of $ \partial \Gamma $.
We aim to show that if $ \Lambda G $ is indivisible in
$ \partial \Gamma $ at each point of $ \Lambda_0 G $, then it
is (globally) indivisible in $ \partial \Gamma $ (Corollary 8.8).
Moreover, if $ \Lambda G(v) $ is indivisible in $ \partial \Gamma $
at some point $ p \in \Lambda G(v) $, then $ \Lambda G $ is
also indivisible in $ \partial \Gamma $ at $ p $ (Proposition 8.9).
As a corollary, we deduce (Corollary 8.10) that if 
$ \Lambda G(v) $ is indivisible in $ \partial \Gamma $
for all $ v \in V(T) $, then $ \Lambda G $
is indivisible in $ \partial \Gamma $.
(Note that this is the essential ingredient in showing that
$ \partial \Gamma $ has no global cut point, as in Lemma 8.5.)
\par
Recall the notation $ \Xi $, $ J(\xi) $, $ H(\xi) $, $ B({\vec e}) $ etc\ from
Section 6.
We begin with the following observation:
\gap
\lemma{8.6}
{\sl
$ \Lambda G $ is indivisible in $ \partial \Gamma $ over
$ \partial \Gamma \bksl \Lambda G $.
}
\gap
\prf
Suppose $ p \in \partial \Gamma \bksl \Lambda G $.
Then, by Lemma 6.1, $ p \in J(\xi) \bksl \fr J(\xi) $ for some
$ \xi \in \Xi $.
Let $ K $ be the closure of $ \partial \Gamma \bksl J(\xi) $ in
$ \partial \Gamma $.
By Lemma 7.2, $ K $ is connected.
Moreover $ \Lambda G \subseteq K $.
Suppose $ O,U \subseteq M $ with $ OpU $.
Without loss of generality, we can suppose that
$ K \cap U = \emptyset $.
(Otherwise $ O \cap K $ and $ U \cap K $ would partition $ K $.)
But $ \Lambda G \subseteq K $, and so $ \Lambda G \cap U = \emptyset $.
\endprf
\par
Recall the notation $ {\script S}_0({\vec e}) $, $ \mathord{\simeq}_S $
etc\ from Section 7.
\par
For each $ {\vec e} \in {\vec E}(T) $, we shall choose
$ S({\vec e}) \in {\script S}_0({\vec e}) $.
We do this equivariantly with respect to the action of
$ G $.
Thus, $ N = \max \{ \diam S({\vec e}) \mid {\vec e} \in {\vec E}(T) \}
< \infty $ (where $ \diam $ denotes diameter with respect to
combinatorial distance in $ T $).
\gap
\lemma{8.7}
{\sl
$ \Lambda G $ is indivisible in $ \partial \Gamma $ over
$ \Lambda_\infty G $.
}
\gap
\prf
Clearly, we can assume that $ \Lambda_\infty G $ is non-empty, and
hence dense in $ \Lambda G $.
Suppose that $ p \in \Lambda_\infty G $, and $ O,U \subseteq
\partial \Gamma $ with $ OpU $.
If $ O \cap \Lambda G \ne \emptyset $, then $ O \cap \Lambda_\infty G
\ne \emptyset $, and similarly for $ U $.
Thus, suppose, for contradiction, that there exist
$ x \in O \cap \Lambda_\infty G $ and $ y \in U \cap \Lambda_\infty G $.
Clearly $ x $, $ y $ and $ p $ are all distinct.
\par
Now, let $ v \in V(T) $ be the median of the points $ x,y,p \in
\partial T $.
In other words, $ v $ is the unique intersection point of the
three arcs connecting the points $ x $, $ y $ and $ p $ pairwise.
Let $ \alpha $ be the ray from $ v $ to $ p $, and let $ w \in V(T) $
be that vertex at distance $ N+1 $ from $ v $ along $ \alpha $.
Let $ {\vec e} $ be the directed edge of $ \alpha $ pointing
towards $ p $ with $ \rhead({\vec e}) = w $ (so that
$ \dist(v,\rtail({\vec e})) = N $.
Thus $ x,y \in \partial \Phi({\vec e}) $ and $ p \in \partial
\Phi(-{\vec e}) $.
\par
Write $ S = S({\vec e}) $, so that $ \diam S \le N < \dist(v,w) $.
Now $ v $ is the nearest point to $ w $ in the biinfinite
arc connecting $ x $ to $ y $.
We see that this arc does not meet $ S $, and so
$ \pi_S x = \pi_S y $.
In particular, $ x \simeq_S y $, and so, since $ S \in
{\script S}_0({\vec e}) $, we have $ x \simeq y $.
By Lemma 7.6, $ x $ and $ y $ lie in the same component of
$ B({\vec e}) $.
But, $ \partial \Phi(-{\vec e}) \cap B({\vec e}) = \emptyset $, and
so $ p \notin B({\vec e}) $.
But this contradicts the fact that $ p $ separates $ x $ from $ y $.
(More formally, $ O \cap B({\vec e}) $ and $ U \cap B({\vec e}) $
partition $ B({\vec e}) $ into two non-empty open sets.)
\endprf
\par
Putting Lemma 8.7 together with Lemma 8.6, we obtain:
\gap
\corol{8.8}
{\sl
If $ \Lambda G $ is indivisible in $ \partial \Gamma $ over
$ \Lambda_0 G $, then $ \Lambda G $ is (globally) indivisible in
$ \partial \Gamma $.
}
\endprf
\par
Next, we show:
\gap
\propn{8.9}
{\sl
If $ \Lambda G(v) $ is indivisible in $ \partial \Gamma $ at the
point $ p \in \Lambda G(v) $, then $ \Lambda G $ is indivisible in
$ \partial \Gamma $ at $ p $.
}
\gap
\prf
First, note that if $ T $ is trivial, then $ G=G(v) $, so there is
nothing to prove.
We can thus assume that $ T $ is non-trivial.
\par
Suppose that $ O,U \subseteq \partial \Gamma $ with $ OpU $.
Since $ \Lambda G(v) $ is indivisible in $ \partial \Gamma $ at $ p $,
we can assume that $ U \cap \Lambda G(v) = \emptyset $.
We claim that $ U \cap \Lambda G = \emptyset $.
Since $ \Lambda_\infty G $ is dense in $ \Lambda G $, it's
enough to show that $ U \cap \Lambda_\infty G = \emptyset $.
\par
Suppose, to the contrary, that there is some $ x \in U \cap
\Lambda_\infty G $.
Let $ G_0 \subseteq G(v) $ be the setwise stabiliser of the
$ \mathord{\sim} $--class of $ x $.
By Lemma 7.7, $ G_0 $ is infinite.
Now a hyperbolic group cannot contain an infinite torsion subgroup
(see for example [\GhH]) and so we can find some $ g \in G_0 $ of
infinite order.
\par
Now, for each $ i \in {\bbf Z} $, $ g^i x \sim x $, so, by
Lemma 7.8, there is a connected subset (in fact a subcontinuum), $ K $,
containing $ x $ and $ g^i x $, with $ K \cap \Lambda G(v) = \emptyset $.
Since $ p \in \Lambda G(v) $, we have $ K \subseteq \partial \Gamma
\bksl \{ p \} $.
Thus, $ K \subseteq U $.
(Otherwise $ O \cap K $ and $ U \cap K $ would partition $ K $.)
In particular, $ g^i x \in U $.
Now, as $ i \rightarrow \infty $, the sequences $ g^i x $ and
$ g^{-i} x $ converge on distinct points, $ a,b \in
\Lambda G_0 \subseteq \Lambda G(v) $.
Since $ U \cup \{ p \} $ is closed, we have $ a,b \in U \cup \{ p \} $,
and so, without loss of generality, $ a \in U $.
But now, $ a \in U \cap \Lambda G(v) $, contradicting the
assumption that $ U \cap \Lambda G(v) = \emptyset $.
\endprf
\par
Putting Proposition 8.9 together with
Corollary 8.8, we get:
\gap
\corol{8.10}
{\sl
Suppose that, for all $ v \in V_{\rm inf}(T) $, the continuum
$ \Lambda G(v) $ is indivisible in $ \partial \Gamma $ over
$ \Lambda G(v) $.
Then, $ \Lambda G $ is (globally) indivisible in $ \partial \Gamma $.
}
\endprf
\par
Of course, it's enough to suppose that each continuum
$ \Lambda G(v) $ has no global cut point.
\par
Finally, putting Corollary 8.10 together with Lemma 8.5, we get
the main result of this section, namely Theorem 8.1.
\gap\gap
\section{Strongly accessible groups}
\gap\par
In this final section, we look once more at the property of strong
accessibility over finite and two-ended subgroups.
We begin with general groups, and specialise to finitely presented
groups.
We finish by showing how Theorem 8.1, together with the results of
[\Boa,\Boc] imply that the boundary of a one-ended strongly accessible
hyperbolic group has no global cut point (Theorem 9.3).
\par
As discussed in the introduction, the issue of strong accessibility
is concerned with sequences of splittings over a class of subgroups
(in particular, the class of finite and two-ended subgroups), and
when such sequences must terminate.
In general, this may depend on the choices of splittings that we make at
each stage of the process.
We first describe a few general results which imply, at least for
finitely presented groups, that we can assume that at any given stage,
we can split over finite groups whenever this is possible.
\par
Suppose, for the moment, that $ \Gamma $ is any group, and that
$ G_1 $ and $ G_2 $ are one-ended subgroups with $ G_1 \cap G_2 $
infinite.
Then the group, $ \bet{G_1 \cup G_2} $, generated by $ G_1 $ and $ G_2 $
is also one-ended.
(For if not, there is a non-trivial action of $ \bet{G_1 \cup G_2} $
on a tree, $ T $, with finite edge stabilisers.
Now, since the groups, $ G_i $ are one-ended, they each fix a unique vertex
of $ T $.
Since $ G_1 \cap G_2 $ is infinite, this must be the same vertex, contradicting
the non-triviality of the action.)
Note that essentially the same argument works if $ G_1 $ is one-ended and
$ G_2 $ is two-ended.
\par
Similarly, suppose that $ G \le \Gamma $ is one-ended, and $ g \in \Gamma $
with $ G \cap gGg^{-1} $ infinite.
Then $ \bet{G,g} $ is one-ended.
(Since if $ \bet{G,g} $ acts on a tree, $ T $, with finite edge stabilisers,
then $ G $ and $ gGg^{-1} $ must fix the same unique vertex of $ T $.
Thus, $ g $ must also fix this vertex, again showing that the action is
trivial.)
Recall that the commensurator, $ \Comm(G) $, of $ G $ is the set of elements
$ g \in \Gamma $ such that $ G \cap gGg^{-1} $ has finite index in $ G $.
Thus, $ \Comm(G) $ is a subgroup of $ \Gamma $ containing $ G $.
We see that if $ G $ is one-ended, then so is $ \Comm(G) $.
\par 
Now, suppose that $ \Gamma $ is accessible over finite groups.
Then every one-ended subgroup of $ \Gamma $ is contained in a unique
maximal one-ended subgroup of $ \Gamma $.
Each maximal one-ended subgroup is equal to its commensurator, and
there are only finitely many conjugacy classes of such subgroups.
If $ G $ is a maximal one-ended subgroup, and $ H \le G $ is two-ended,
then either $ H \le G $ or else $ H \cap G $ is finite.
Moreover, $ H $ can lie in at most one maximal one-ended subgroup.
These observations follow from the remarks of the previous two paragraphs.
They can also be deduced by considering the action of $ H $ on
a complete $ \Gamma $--tree.
\par
Now, suppose that $ \Gamma $ splits as an amalgamated free product or
HNN--extension over a two-ended subgroup.
This corresponds to a $ \Gamma $--tree, $ \Sigma $, with just one orbit
of edges, and with two-ended edge stabiliser.
We consider two cases, depending on whether or not the edge group is
elliptic or hyperbolic, ie\ whether or not it lies in a one-ended subgroup
of $ \Gamma $.
\par
Consider, first, the case where the edge stabiliser of $ \Sigma $ does
not lie in a one-ended subgroup, and hence intersects every one-ended
subgroup in a finite group.
In this case, we have:
\gap
\lemma{9.1}
{\sl
Suppose $ v \in V(\Sigma) $.
Then, each maximal one-ended subgroup of $ \Gamma(v) = \Gamma_\Sigma(v) $
is a maximal one-ended subgroup of $ \Gamma $.
Moreover, every maximal one-ended subgroup of $ \Gamma $ arises in
this way (for some $ v \in V(\Sigma) $).
}
\gap
\prf
Suppose, first, that $ G $ is any one-ended subgroup of $ \Gamma $.
Then, $ G $ must lie inside some (unique) vertex stabiliser
$ \Gamma(v) $.
(Otherwise, $ G $ would split over a group of the form $ G \cap H $,
where $ H $ is an edge-stabiliser.
But $ G \cap H $ is finite, contradicting the fact that $ G $ is
one-ended.)
If $ G $ is maximal in $ \Gamma $, then clearly it is also
maximal in $ \Gamma(v) $.
\par
Conversely, suppose that $ G $ is a maximal one-ended
subgroup of a vertex stabiliser, $ \Gamma(v) $.
Let $ G' $ be the unique maximal one-ended subgroup of $ \Gamma $
containing $ G $.
By the first paragraph, $ G' $ lies inside some vertex group, which
must, in this case, be $ \Gamma(v) $.
By maximality in $ \Gamma(v) $, we must therefore have $ G=G' $.\Endprf
\par
The second case is when an edge group lies inside some one-ended
subgroup.
To consider this case, fix and edge $ e $ of $ \Sigma $, with
endpoints $ v,w \in V(\Sigma) $.
Now, $ \Gamma(e) $ lies inside a unique maximal one-ended subgroup,
$ \Gamma_0 $, of $ \Gamma $.
Any other maximal one-ended subgroup of $ \Gamma $ must
intersect $ \Gamma(e) $ in a finite subgroup.
In this case, we have:
\gap
\lemma{9.2}
{\sl
$ \Gamma_0 $ splits as an amalgamated free product or HNN extension
over $ \Gamma(e) $, with incident vertex groups equal to
$ \Gamma_0 \cap \Gamma(v) $ and $ \Gamma_0 \cap \Gamma(w) $.
Each maximal one-ended subgroup of $ \Gamma(v) $ is a maximal one-ended
subgroup of $ \Gamma_0 \cap \Gamma(v) $ or of $ \Gamma $
(and similarly for $ w $).
Every maximal one-ended subgroup of $ \Gamma_0 \cap \Gamma(v) $
arises in this way.
Each maximal one-ended subgroup of $ \Gamma $ is conjugate, in $ \Gamma $,
to $ \Gamma_0 $ or to a maximal one-ended subgroup of $ \Gamma(v) $
or $ \Gamma(w) $.
}
\gap
\prf
Suppose $ G $ is a maximal one-ended subgroup of $ \Gamma $.
Either $ G $ contains some edge-stabiliser, so that some conjugate
of $ G $ contains $ \Gamma(e) $ and hence equals $ \Gamma_0 $,
or else $ G $ meets each edge stabiliser in a finite group.
In the latter case, we see, as in Lemma 9.1, that $ G $ is a
maximal one-ended subgroup of a vertex group.
\par
Now suppose that $ G $ is a maximal one-ended subgroup of
$ \Gamma(v) $.
Let $ G' $ be the maximal one-ended subgroup of $ \Gamma $
containing $ G $.
 From the first paragraph, we see that either $ G' = \Gamma_0 $,
or $ G' $ is a maximal one-ended subgroup of $ \Gamma(v) $.
In the former case, we see that $ G \subseteq \Gamma_0 \cap \Gamma(v) $,
and must therefore be maximal one-ended in $ \Gamma_0 \cap \Gamma(v) $.
The latter case, we obtain $ G = G' $.
\par
Finally suppose that $ G $ is a maximal one-ended subgroup of
$ \Gamma_0 \cap \Gamma(v) $.
Let $ G' $ be the maximal one-ended subgroup of $ \Gamma(v) $
containing $ G $.
 From the previous paragraph, we see that
$ G' \subseteq \Gamma_0 \cap \Gamma(v) $, so $ G = G' $.
\par
It remains to show that $ \Gamma_0 $ splits over $ \Gamma(e) $
in the manner described.
This amounts to showing that if $ H $ is an edge stabiliser and
a subgroup of $ \Gamma_0 \cap \Gamma(v) $, then $ H $ is
conjugate in $ \Gamma_0 \cap \Gamma(v) $ to $ \Gamma(e) $,
(and similarly for $ w $).
\par
We know that there must be some $ g \in \Gamma(v) $ such
that $ H = g \Gamma(e) g^{-1} $
Now, $ H \subseteq \Gamma_0 \cap g \Gamma_0 g^{-1} $.
Since $ H $ is infinite, it follows that the group generated
by $ \Gamma_0 $ and $ g \Gamma_0 g^{-1} $ must be one-ended, and so,
by maximality, must equal $ \Gamma_0 $.
Hence, $ g \Gamma_0 g^{-1} = \Gamma_0 $.
In particular, $ g \in \Comm(\Gamma_0) $.
But, from the earlier discussion, $ \Comm(\Gamma_0) = \Gamma_0 $,
and so $ g \in \Gamma_0 \cap \Gamma(v) $ as required.
\endprf
\par
We now go on to describe the notion of strong accessibility.
To set up the notation, let $ \Gamma $ be any group, and let
$ {\script C} $ be any conjugacy--invariant set of subgroups of $ \Gamma $.
(In the case of interest, $ {\script C} $ will be the set of all finite
and two-ended subgroups of $ \Gamma $.)
We want to look at sequences of splittings of $ \Gamma $ over $ {\script C} $,
where the only information retained at each stage will be the vertex groups of
the previous splittings.
In other words, we get a sequence of conjugacy invariant sets of subgroups
of $ \Gamma $.
(In fact, if $ {\script C} $ is closed under isomorphism, we can just
view these as isomorphism classes of groups.)
Note that finite groups can never split non-trivially, and so for
our purposes, we can throw away finite subgroups whenever they arise.
\par
To be more formal, suppose that $ {\script J} $ and $ {\script J}' $ are
both conjugacy invariant sets of subgroups of $ \Gamma $.
We say that $ {\script J}' $ is obtained by splitting $ {\script J} $ over
$ {\script C} $ if it has the form
$ {\script J}' = \bigcup_J {\script J}(J) $, where $ {\script J}(J) $
is the set of ($ \Gamma $--conjugacy classes of) infinite vertex groups
of some splitting of $ J $ as a finite graph of groups over $ {\script C} $,
and where $ J $ ranges over a conjugacy transversal in $ {\script J} $.
Thus, a sequence of splittings of $ \Gamma $ over $ {\script C} $ consists
of a sequence, $ {\script J}_0, {\script J}_1, {\script J}_2, \ldots $,
where $ {\script J}_0 = \{ \Gamma \} $, and each $ {\script J}_{i+1} $
is obtained as a splitting of $ {\script J}_i $ over $ {\script C} $ in
the manner just described.
Note that, by induction, each of the sets $ {\script J}_i $ is a finite
union of conjugacy classes in $ \Gamma $.
Note also that we can assume, if we wish, by introducing some intermediate
steps, that each $ {\script J}_{i+1} $ is obtained from $ {\script J}_i $ by
splitting one of the conjugacy classes of $ {\script J}_i $ as
an amalgamated free product or HNN extension, while leaving the
remaining groups unchanged.
We say that the sequence terminates, if for some $ n $, none of
the elements of $ {\script J}_n $ split non-trivially over $ {\script C} $.
We say that $ \Gamma $ is \em{strongly accessible} over $ {\script C} $ if
there exists such a sequence which terminates.
\par
Suppose that $ {\script J} $ is a union of conjugacy classes of
subgroups of $ \Gamma $, each accessible over finite groups.
Let $ {\script F}({\script J}) = \bigcup_{J \in {\script J}} {\script F}(J) $,
where $ {\script F}(J) $ is the set of maximal one-ended subgroups of $ J $.
Thus $ {\script F}({\script J}) $ is obtained by $ {\script J} $ by
splitting over the class of finite subgroups of $ \Gamma $, in the
sense defined above.
\par
Let us now suppose that $ \Gamma $ is finitely presented, and that
$ {\script C} $ is the set of all finite and one-ended subgroups of
$ \Gamma $.
Suppose that $ ({\script J}_i)_i $ is a sequence of splitting of $ \Gamma $
over $ {\script C} $.
By induction, each element of each $ {\script J}_i $ is finitely presented
and hence accessible over finite groups.
We can thus form a sequence $ ({\script F}_i)_i $ where $ {\script F}_i
= {\script F}({\script J}_i) $.
Now, we can assume that $ {\script J}_{i+1} $ is obtained from
$ {\script J}_i $ by splitting an element of $ {\script J}_i $ as
an amalgamated free product or HNN extension either over a finite group
or over a two-ended group.
In the former case, we see that $ {\script F}_{i+1} = {\script F}_i $.
In the latter case, we see, from Lemmas 9.1 and 9.2, that $ {\script F}_{i+1} $
is obtained from $ {\script F}_i $ by first splitting some element over
a two-ended subgroup, and then, if necessary splitting over some finite
subgroups to reduce ourselves again to one-ended groups.
Thus, after inserting some intermediate steps if necessary, we can suppose
that the sequence $ ({\script F}_i)_i $ is also a sequence of splittings
of $ \Gamma $ over $ {\script C} $.
If the sequence $ ({\script J}_i)_i $ terminates at $ {\script J}_n $, then
$ {\script F}_n = {\script F}({\script J}_n) = {\script J}_n $, so
$ ({\script F}_i)_i $ also terminates (and in the same set of subgroups).
\par
In summary, we see that if $ \Gamma $ is finitely presented, and strongly
accessible over $ {\script C} $, then we can find a terminating sequence
of splittings over $ {\script C} $ where we split over finite groups
wherever possible (in priority to splitting over two-ended subgroups).
In other words, we only ever need to split one-ended groups over two-ended
subgroups and to split infinite-ended and two-ended groups over finite
subgroups.
\par
Finally, suppose that $ \Gamma $ is a strongly accessible one-ended hyperbolic
group, and that $ {\script J}_0, {\script J}_1 ,\ldots, {\script J}_n $
is a sequence of splitting of $ \Gamma $ over finite and one-ended
subgroups, which terminates in $ {\script J}_n $.
In this case, each elements of each $ {\script J}_i $ is quasiconvex, and
hence intrinsically hyperbolic.
Moreover, we can suppose, as above, that the only groups we ever split over
two-ended groups are one-ended.
\par
Now, each element of $ {\script J}_n $ is one-ended and does not split
over any two-ended subgroup.
 From the results of [\Boa,\Boc], we see that each element of $ {\script J}_n $
has no global cut point in its boundary.
Now, applying Theorem 8.1 inductively, we conclude that this is also
true of $ \Gamma $.
\par
We have shown:
\gap
\theorem{9.3}
{\sl
Suppose that $ \Gamma $ is a one-ended hyperbolic group which is strongly
accessible over finite and two-ended subgroups.
Then, $ \partial \Gamma $ has no global cut point.
}
\endprf
\par
As mentioned in the introduction, Delzant and Potyagailo have shown that
every finitely presented group, $ \Gamma $, is strongly accessible over
any ``elementary'' class of subgroups, $ {\script C} $.
In particular, this deals with the case where $ \Gamma $ is hyperbolic, and
where $ {\script C} $ is the set of finite and two-ended subgroups of
$ \Gamma $.
We thus conclude that the boundary of any one-ended hyperbolic group has
no global cut point, and is thus locally connected by the result of [\BeM].

\np
\references

\ppar

\bigskip
{\small \parskip 0pt \leftskip 0pt \rightskip 0pt plus 1fil \def\\{\par}
\sl\theaddress\par
\medskip
\rm Email:\stdspace\tt\theemail\par}

{\small Received:\qua 15 November 1997\qquad Revised:\qua 10 August 1998}

\end